\documentclass[envcountsame,envcountsect]{svmult}

\usepackage{makeidx}         
\usepackage{graphicx}        
\usepackage{multicol}        
\usepackage[bottom]{footmisc}

\makeindex             
\usepackage{amsmath}
\usepackage{amssymb}
\newtheorem{Number}[theorem]{\!\!}
\newenvironment{numba}{\begin{Number}\rm}{\end{Number}}
\newcommand{\wb}{\overline}
\newcommand{\ve}{\varepsilon}
\newcommand{\wt}{\widetilde}

\newcommand{\isom}{\cong}
\newcommand{\mto}{\mapsto}
\newcommand{\N}{{\mathbb N}}
\newcommand{\bS}{{\mathbb S}}
\newcommand{\A}{{\mathbb A}}
\newcommand{\R}{{\mathbb R}}
\newcommand{\C}{{\mathbb C}}
\newcommand{\K}{{\mathbb K}}
\newcommand{\Q}{{\mathbb Q}}

\newcommand{\TOP}{{\mathbb T}{\mathbb O}{\mathbb P}}
\newcommand{\LIE}{{\mathbb L}{\mathbb I}{\mathbb E}}
\newcommand{\cV}{{\mathcal V}}
\newcommand{\cO}{{\mathcal O}}
\newcommand{\cS}{{\mathcal S}}

\newcommand{\cT}{{\mathcal T}}
\newcommand{\cH}{{\mathcal H}}
\newcommand{\gl}{{\mathfrak{g}\mathfrak{l}}}
\newcommand{\cg}{{\mathfrak g}}

\DeclareMathOperator{\Hom}{Hom}
\newcommand{\wh}{\widehat}
\newcommand{\sub}{\subseteq}
\DeclareMathOperator{\GL}{GL}

\DeclareMathOperator{\id}{id}

\DeclareMathOperator{\Diff}{Diff}
\DeclareMathOperator{\GermDiff}{GermDiff}
\DeclareMathOperator{\Supp}{supp}

\DeclareMathOperator{\Hol}{Hol}
\DeclareMathOperator{\SO}{SO}
\DeclareMathOperator{\Sp}{Sp}
\DeclareMathOperator{\SU}{SU}
\DeclareMathOperator{\U}{U}
\DeclareMathOperator{\gO}{O}
\DeclareMathOperator{\Germ}{Germ}
\DeclareMathOperator{\evol}{evol}
\DeclareMathOperator{\Lie}{{\bf L}}
\DeclareMathOperator{\lcx}{lcx}
\DeclareMathOperator{\cts}{cts}
\DeclareMathOperator{\DL}{DL}
\DeclareMathOperator{\DLG}{DLG}
\DeclareMathOperator{\Gau}{Gau}
\DeclareMathOperator{\conv}{conv}
\newcommand{\SET}{\mbox{${\mathbb S}{\mathbb E}{\mathbb T}$}}
\newcommand{\G}{\mathbb{G}}
\newcommand{\dl}{{\displaystyle\lim_{\longrightarrow}}}
\newcommand{\pl}{{\displaystyle\lim_{\longleftarrow}}}
\newcommand{\semid}{\mbox{$\times\!$\rule{.15 mm}{1.83 mm}}}
\begin{document}

\title*{Direct limits of infinite-dimensional Lie groups}
\author{Helge Gl\"{o}ckner\inst{1}}
\institute{Universit\"{a}t Paderborn, Institut f\"{u}r
Mathematik, Warburger Str.\ 100,\\
33098 Paderborn, Germany.
\texttt{glockner@math.uni-paderborn.de}}
\maketitle
\begin{abstract}
Many
infinite-dimensional
Lie groups~$G$ of interest
can be expressed as the union $G=\bigcup_{n\in \N}G_n$
of an ascending sequence
$G_1\sub G_2\sub\cdots$ of (finite- or infinite-dimensional)
Lie groups. In this survey article, we compile general
results\linebreak
concerning
such ascending unions,
describe the main classes of examples,
and\linebreak
explain what the general theory tells us about these.\\[2mm]
Classification. 22E65
(Primary) 
26E15; 
46A13; 
46G20; 
46T05; 
46T20; 
46T25; 
54B35; 
54D50; 
55Q10; 
58B05; 
58D05 (Secondary)\\[2mm] 
\end{abstract}
\section{\,Introduction}
Many infinite-dimensional
Lie groups $G$ can be expressed as the union
$G=\bigcup_{n\in \N}\,G_n$
of a sequence $G_1\sub G_2\sub \cdots$
of (finite- or infinite-dimensional) Lie groups,
such that the inclusion maps $j_n\colon G_n\to G$
and $j_{m,n}\colon G_n\to G_m$ (for $n\leq m$)
are smooth homomorphisms.
Typically, the steps $G_n$
are Lie groups of a simpler type,
and one hopes (and often succeeds)
to deduce results concerning~$G$
from information available
for the Lie groups~$G_n$.\\[2.5mm]
The goals of this article are twofold:
\begin{itemize}
\item
To survey general results
on ascending unions of Lie groups
and their properties;\vspace{0.5mm}
\item
To collect concrete classes of examples
and explain how the general theory
specializes in these cases.
\end{itemize}
One
typical class of examples
is given by the groups $\Diff_c(M)$\index{diffeomorphism group}
of smooth diffeomorphisms
$\phi\colon M\to M$
of $\sigma$-compact, finite-dimensional smooth manifolds~$M$
which are compactly supported
in the sense that the set
$\{x\in M\colon \phi(x)\not= x\}$
has compact
closure. The group operation is
composition of diffeomorphisms.
It is known that $\Diff_c(M)$
is a Lie group (see \cite{Mic} or \cite{DIF}).
Furthermore,
\[
\Diff_c(M)\;=\: \bigcup_{n\in \N} \, \Diff_{K_n}(M)\vspace{-2mm}
\]
for each exhaustion
$K_1\sub K_2\sub\cdots$
of~$M$ by compact sets
(with $K_n$ in the interior of $K_{n+1}$),
where $\Diff_{K_n}(M)$ is the Lie group
of smooth diffeomorphisms of~$M$ supported in~$K_n$.
The manifold structure of $\Diff_c(M)$
is modelled on the space
$\cV_c(M)$ of compactly supported
smooth vector fields,
which is an LF-space with a
complicated topology.
By contrast, $\Diff_{K_n}(M)$
is modelled on the space
$\cV_{K_n}(M)$ of smooth vector fields
supported in~$K_n$,
which is a Fr\'{e}chet space.
Many specific tools of infinite-dimensional
calculus can be applied to $\cV_{K_n}(M)$,
e.g.\ to clarify differentiability
questions for functions on this space.

In other typical cases,
each $G_n$ is finite-dimensional
(a particularly well-understood situation)
or modelled on a Banach space,
whence again special tools
are available to deal with the Lie groups
$G_n$ (but not \emph{a priori} for~$G$).

Besides diffeomorphism groups,
we shall also discuss the following
major classes of examples
(described in more detail in Section~\ref{secuninf}):\vspace{.8mm}
\begin{itemize}
\item
The ``test function groups''\index{test function group}
$C^\infty_c(M,H)=\bigcup_{n\in \N}\,C_{K_n}^\infty(M,H)$
of compactly supported Lie group-valued
smooth mappings on a $\sigma$-compact
smooth manifold $M=\bigcup_{n\in \N}K_n$;\vspace{.5mm}
\item
Weak direct products $\prod_{n\in \N}^* H_n:=\bigcup_{n\in \N}
\prod_{k=1}^n H_k$\index{weak direct product}
of Lie groups $H_n$;\vspace{.5mm}
\item
Unions $A^\times=\bigcup_{n\in \N}A_n^\times$
of unit groups
of Banach algebras $A_1\sub A_2\sub\cdots$;\vspace{.5mm}
\item
The groups $\Germ(K,H)$ of germs\label{group of germs of analytic maps}
of analytic mappings on open neighbourhoods\index{Sobolev--Lie group}
of a compact subset $K$ of a metrizable
complex locally convex space,
with values in a complex Banach--Lie group~$H$.
\item
The group
$H^{\downarrow s}(K,F)
=\bigcup_{t>s}H^t(K,F)=
\bigcup_{n\in \N}H^{s+\frac{1}{n}}(K,F)$,
where $K$ is a compact smooth manifold,
$s\geq \dim(K)/2$,
$F$ a finite-dimensional Lie group, and
$H^t(K,F)\sub C(K,F)$ the integral subgroup
whose Lie algebra is the
Sobolev space
$H^t(K,\Lie(F))$
of functions with values in the Lie algebra
$\Lie(F)$ of~$F$.\vspace{.8mm}
\end{itemize}

For $s=\dim(K)/2$, the Lie group $H^{\downarrow s}(K,F)$
is particularly interesting,
because a Hilbert--Lie group
$H^s(K,F)$ is not available in this case.
In some situations,
$H^{\downarrow s}(K,F)$ may serve as a substitute
for the missing group.

We shall also discuss the group $\GermDiff(K,X)$\index{group of germs
of diffeomorphisms}
of germs of analytic diffeo\-morphisms $\gamma$
around a compact set $K$ in a finite-dimensional
complex vector space~$X$,\vspace{-.3mm}
such that $\gamma|_K=\id_K$.
This group is not considered as a union of groups, but
as a union
of Banach \emph{manifolds} $M_1\sub M_2\sub\cdots$.

Among others, we shall discuss
the following topics
in our general setting
(and for the preceding examples):\vspace{.8mm}
\begin{itemize}
\item
Direct limit properties
of ascending unions;\vspace{.5mm}
\item
Homotopy groups
of ascending unions;\vspace{.5mm}
\item
When ascending unions
are regular Lie groups in
Milnor's sense;\vspace{.5mm}
\item
Questions concerning subgroups of ascending unions.\vspace{.8mm}
\end{itemize}

We now describe the main problems
and questions
in more detail, together
with some essential concepts.
As a rule, references to the literature,
answers (and partial
answers) will only be given later,
in the actual article.
\subsection{Direct limit properties of ascending unions}\label{dlpropannounce}
Consider a Lie group~$G$ which is an
ascending union $G=\bigcup_{n\in \N}\, G_n$
of Lie groups, and a map $f\colon G\to X$.
It is natural to ask:\index{direct limit properties}
\begin{description}
\item[(a)]
If~$X$ is a smooth manifold (modelled on a locally convex
space) and $f|_{G_n}$ is smooth
for each $n\in \N$, does it follow that $f$
is smooth?
\item[(b)]
If~$X$ is a topological space and $f|_{G_n}$ is continuous
for each $n\in \N$, does it follow that $f$
is continuous?
\item[(c)]
If~$X$ is a Lie group, $f$ is a homomorphism
of groups and $f|_{G_n}$ is smooth
for each $n\in \N$, does it follow that $f$
is smooth?
\item[(d)]
If~$X$ is a topological group,
$f$ is a homomorphism and $f|_{G_n}$ is continuous
for each $n\in \N$, does it follow that $f$
is continuous?
\end{description}

As we shall see, (a) and (b) are frequently not true
(unless compactness can be brought into play),
while (c) and (d) hold for our typical examples.

The preceding questions can be re-cast
in category-theoretic terms:
They amount to asking if $G$ is the direct limit
$\dl\,G_n$\vspace{-.6mm}
in the categories of smooth manifolds,
topological spaces, Lie groups,
resp.,
topological groups
(see~\ref{explainnow}).
The relevant concepts from
category theory
will be recalled in Section~\ref{secprel}.

Questions (b) and (d) can be asked
just
as well if $G$ and each $G_n$ merely is a topological
group, and each inclusion map is a continuous
homomorphism.
Essential progress concerning
direct limits of topological groups
and their\linebreak
relations to direct limits
of topological spaces
were achieved in the last ten years,
notably by N. Tatsuuma, E. Hirai, T. Hirai
and N. Shimomura (see \cite{TSH} and \cite{HSTH})
as well as A. Yamasaki~\cite{Yam}.
In Section~\ref{secgp},
we recall the most relevant results.
\subsection{Existence of direct limit charts --
an essential hypothesis}\label{susecdlcha}
Meaningful results concerning the topics
raised above can only be expected
under additional hypotheses.
For instance,
our general setting includes
the situation
where each $G_n$ is discrete
but $G$ is not
(as we only assume that the inclusion maps
$G_n\to G$ are smooth).
In this situation,
algebraic properties
of the groups $G_n$
(like simplicity or perfectness)
pass to~$G$,
but we cannot expect to gain
information concerning the topological or
differentiable structure of~$G$
from information on
the groups~$G_n$.

A very mild additional hypothesis
is the existence
of a \emph{direct limit chart}.\index{direct limit chart}
Roughly speaking, this is
a chart of $G$ juxtaposed from charts of the Lie groups~$G_n$.
The formal definition
reads as follows (cf.\ \cite[Definition 2.1]{COM}):\\[2.7mm]
{\bf Definition.} \,A Lie group $G=\bigcup_{n\in \N}G_n$ is said to
admit a \emph{weak direct limit chart}\index{weak direct limit chart}
if there exists $n_0\in \N$,
charts $\phi_n\colon U_n\to V_n$
from open identity neighbourhoods
$U_n\sub G_n$ onto open $0$-neighbourhoods
$V_n\sub \Lie(G_n)$ in the tangent space
$\Lie(G_n):=T_{\bf 1}(G_n)$ at~${\bf 1}$ for $n\geq n_0$
and a chart $\phi\colon U\to V$
from an open\linebreak
identity neighbourhood
$U\sub G$ onto an open $0$-neighbourhood
$V\sub \Lie(G)$,
such that
\begin{description}
\item[(a)]
$U=\bigcup_{n\geq n_0}U_n$
and $U_n\sub U_{n+1}$ for each integer $n\geq n_0$; and
\item[(b)]
$\phi_{n+1}|_{U_n}=\Lie(j_{n+1,n})\circ \phi_n$
and
$\phi|_{U_n}=\Lie(j_n)\circ \phi_n$
for each $n\geq n_0$.
\end{description}
If, furthermore,
$\Lie(G)=\dl\,\Lie(G_n)$\vspace{-.5mm}
as a locally convex space,\footnote{Here, we use the bonding maps
$\Lie(j_{m,n})\colon \Lie(G_n)\to\Lie(G_m)$
and the limit maps $\Lie(j_n)\colon \Lie(G_n)\to\Lie(G)$.}
then $G=\bigcup_{n\in\N}G_n$ is said to
admit a \emph{direct limit chart}.

Note that (b) implies that the linear maps
$\Lie(j_n)$ and $\Lie(j_{n+1,n})$
are injective on some $0$-neighbourhood and thus
injective. Hence, identifying
$\Lie(G_n)$ with its image
under $\Lie(j_n)$ in $\Lie(G)$,
we can re-write (b) as
\begin{description}
\item[(b)${}'$]
$\phi|_{U_n}=\phi_n$
and $\phi_{n+1}|_{U_n}=\phi_n$,
for each $n\geq n_0$.
\end{description}
Furthermore, we now simply have
$V=\bigcup_{n\geq n_0}V_n$.

To assume the existence of a direct
limit chart is
a natural requirement,
which is satisfied by all of our main
examples.
It provides a link
between the topologies (resp., manifold structures)
on $G$ and the Lie groups~$G_n$,
and will be encountered in connection with
most of the topics from above.\vspace{-1mm}
\subsection{Homotopy groups of ascending unions
of Lie groups}\label{subsechom}
Given a Lie group $G=\bigcup_{n\in\N}G_n$,
it is natural to ask if
its $k$-th homotopy group
can be calculated in terms of\index{homotopy group}
the homotopy groups $\pi_k(G_n)$ in the form\vspace{-.7mm}
\begin{equation}\label{dlhomot}
\pi_k(G)\;=\; \dl\, \pi_k(G_n)\,,\vspace{-2mm}
\end{equation}
for each $k\in \N_0$.
This is quite obvious if
$G=\bigcup_{n\in \N}G_n$ is \emph{compactly regular}\index{compactly regular}
in the sense that each compact subset~$K$ of~$G$\index{compact regularity}
is a compact subset of some~$G_n$
(see \cite[Proposition~3.3]{HOM};
cf.\ \cite[Remark~3.9]{FUN}
and \cite[Lemma~A.7]{NeO}
for special cases, as well as works
on stable homotopy theory and $K$-theory).

There is another, non-trivial condition:
\emph{If $G=\bigcup_{n\in\N}G_n$
admits a weak\linebreak
direct limit chart,
then} (\ref{dlhomot}) \emph{holds} \cite[Theorem~1.2]{HOM}.
A variant of this condition
even applies if $\bigcup_{n\in \N}G_n$
is merely dense in~$G$ (see Theorem~1.13 in \cite{HOM}).
Moreover, ascending unions
can be replaced with directed unions
over uncountable
families, and Lie groups with manifolds
(see Section~\ref{sechomotop}).
These results are based on approximation arguments.
Analogous results for open subsets of locally
convex spaces are classical~\cite{Pa2}.

We mention that knowledge of $\pi_0(G)=G/G_0$,
the fundamental group $\pi_1(G)$
and $\pi_2(G)$ is essential for the extension theory
of~$G$.
It is needed to understand the Lie group
extensions ${\bf 1} \to A\to \wh{G}\to G\to {\bf 1}$
of~$G$ with abelian kernel,
by recent results of K.-H. Neeb
(see \cite{NeC}, \cite{NeA}, and \cite{NeN}).\vspace{-1mm}
\subsection{Regularity in Milnor's sense}\label{ssecreg}
Roughly speaking, a Lie group $G$
(modelled on a locally convex space)
is called a regular Lie group
if all differential equations on $G$
which are of relevance
for Lie theory can be solved,
and their solutions depend
smoothly on parameters.
To make this more precise,
given $g,h \in G$ and $v\in T_h(G)$
let us write $g\cdot v:=(T_h\lambda_g)(v)\in T_{gh}(G)$,
where $\lambda_g\colon G\to G$,
$x\mto gx$ denotes left translation by~$g$.\\[2.5mm]
{\bf Definition.} A Lie group~$G$
modelled on a locally convex space
is called a \emph{regular Lie group} (in Milnor's
sense) if for each smooth curve\index{regular Lie group}
$\gamma\colon [0,1]\to \Lie(G)$,
there exists a (necessarily unique)
smooth curve $\eta= \eta_\gamma\colon [0,1]\to G$
(a so-called ``product integral'')\index{product integral}
which solves the initial value problem
\begin{eqnarray}
\eta(0) & = & {\bf 1} \label{inval1}\\
\eta'(t) & = & \eta(t)\cdot \gamma(t)\;\;\mbox{for all $t\in [0,1]$}\label{inval2}
\end{eqnarray}
(with ${\bf 1} \in G$ the identity element),
and the ``evolution map''\index{evolution map}
\begin{equation}\label{basicevol}
\evol\colon C^\infty([0,1],\Lie(G))\to G\,,\quad
\evol(\gamma):=\eta_\gamma(1)
\end{equation}
is smooth (see \cite{Mr84}, \cite{GaN}
and \cite{NeS}).
Regularity is a useful property,
which provides a link between $G$
and its Lie algebra.
In particular, regularity
ensures the existence of a smooth exponential
map $\exp_G\colon \Lie(G)\to G$, i.e.,\index{exponential map}
a smooth map such that, for each $v\in \Lie(G)$,\index{exponential function}
\[
\gamma_v\colon \R\to G\,,\quad t\mto \exp_G(tv)
\]
is a homomorphism of groups with initial
velocity $\gamma_v'(0)=v$ (cf.\ \cite{Mr84}).

The modelling space $E$ of a regular Lie group
is necessarily \emph{Mackey\linebreak
complete}
in the sense that the Riemann integral\index{Mackey complete}
$\int_0^1 \gamma(t)\,dt$ exists in~$E$
for each smooth curve $\gamma\colon \R\to E$
(cf.\ Lemma~A.5\,(1) and p.\,4 in \cite{NaW}).
Lie groups modelled
on non-Mackey complete locally convex spaces
need not even have an exponential
map. For example, this pathology occurs
for group $G=A^\times$
of invertible elements
in the normed algebra $A\sub C[0,1]$ of all
restrictions to $[0,1]$
of rational functions without poles in $[0,1]$
(with the supremum norm).
Because~$A^\times$ is an open subset of~$A$,
it is a Lie group,
and it is not hard to see
that a smooth homomorphism
$\gamma_v\colon \R\to A^\times$
with $\gamma'(0)=v$
exists for $v\in A=\Lie(A^\times)$
only if $v$ is a constant function
\cite[Proposition~6.1]{ALG}.
Further information concerning
Mackey completeness can be found \cite{KaM}.

At the time of writing, it is unknown
whether non-regular Lie groups modelled
on Mackey complete
locally convex spaces exist.
However, there is no general method of proof;
for each individual class of Lie groups,
very specific arguments are required to
verify regularity.

\noindent
It is natural to look for
conditions ensuring that a union $G=\bigcup_{n\in \N}G_n$
is a regular Lie group if so is each~$G_n$.
Already the
case of finite-dimensional Lie groups $G_n$ is
not easy~\cite{FUN}. In Section~\ref{secregu},
we preview work in progress concerning
the general case. We also describe a construction
which might lead to non-regular Lie groups
(Proposition~\ref{pathrg}).
The potential counterexamples are weak direct
products of suitable regular Lie groups.
\subsection{Subgroups of direct limit groups}
It is natural to try to use
information concerning the subgroups of Lie groups~$G_n$
to deduce results
concerning the subgroups of a Lie group
$G=\bigcup_{n\in \N}G_n$.
Aiming at a typical
example, let us recall that a topological group~$G$
is said to \emph{have no small subgroups}\index{small subgroups}
if there exists an identity neighbourhood
$U\sub G$ containing no subgroup
of~$G$ except for the trivial subgroup.
Although finite-dimensional
(and Banach-) Lie groups do not have
small subgroups, already for Fr\'{e}chet--Lie groups
the situation changes:
The additive group
of the Fr\'{e}chet space $\R^\N$
has small subgroups.
In fact, every
$0$-neighbourhood contains $]{-r},r[^n\times \R^{\{n+1,n+2,\ldots\}}$
for some $n\in \N$ and $r>0$.
It therefore contains the non-trivial
subgroup $\{0\}\times \R^{\{n+1,n+2,\ldots\}}$.

It is natural to ask whether
a Lie group $G=\bigcup_{n\in \N}G_n$
does not have small subgroups
if none of the Lie groups $G_n$
has small subgroups.
In Section~\ref{secsub}, we describe the available
answers to this question,
and various other results
concerning subgroups of direct limit groups.
\subsection{Constructions of Lie group structures
on ascending unions}
So far, we assumed that $G$ is already equipped
with a Lie group structure.
Sometimes, only an ascending
sequence $G_1\sub G_2\sub \cdots$ of Lie groups
is given such that all inclusion maps $G_n\to G_{n+1}$
are smooth homomorphisms.
It is then natural to ask
whether the union
$G=\bigcup_{n\in \N}G_n$
can be given a Lie group structure making each
inclusion map $G_n\to G$ a smooth homomorphism.\footnote{Or even making $G$
the direct limit $\dl\,G_n$\vspace{-.7mm} in the category
of Lie groups.}
We shall also discuss this
complementary problem (in Section~\ref{secconstr}).
If each $G_n$ is finite-dimensional,
then a Lie group structure on $G$ is always
available.
\subsection{Properties of locally convex direct limits}\label{subseclcx}
To enable an understanding of direct limits
of Lie groups, an understanding of various
properties of locally convex direct
limits is essential,
i.e., of direct limits
in the category of locally convex spaces.\index{locally convex direct limit}
For instance, we shall see that if a Lie group
$G=\bigcup_{n\in \N}G_n$ admits a direct limit chart,
then $G=\dl\,G_n$ as a topological space
if and only if $\Lie(G)=\dl\,\Lie(G_n)$\vspace{-.7mm}
as a topological space.
The latter property
is frequently easier to prove (or refute)
than the first.
Also compact regularity of~$G$
(as in~\ref{subsechom} above)
can be checked on the level of
the modelling spaces (see Lemma~\ref{thsreg}).
Another property is useful:
Consider a locally convex space~$E$
which is a union $\bigcup_{n\in \N}E_n$
of locally convex spaces, such that
all inclusion maps are continuous
linear maps. We say that $E$ is \emph{regular}
(or \emph{boundedly regular}, for added clarity)\index{boundedly regular}
if every bounded subset of~$E$
is a bounded subset of some~$E_n$.
If one wants
to prove that a Lie group $\bigcup_nG_n$ is regular
in Milnor's sense,
then it helps a lot if
one knows that $\Lie(G)=\bigcup_n\Lie(G_n)$
is compactly or boundedly regular
(see Section~\ref{secregu}).
\subsection{Further comments, and some
historical remarks}\label{introhis}
The most typical examples of direct limits of finite-dimensional
Lie groups are unions
of classical groups like
$\GL_\infty(\C)=\bigcup_{n\in\N}\GL_n(\C)$
and its subgroups
$\GL_\infty(\R)=\bigcup_{n\in\N}\GL_n(\R)$,
$\gO_\infty(\R)=\bigcup_{n\in \N}\gO_n(\R)$
and $\U_\infty(\C)=\bigcup_{n\in \N}\U_n(\C)$,
where $A\in \GL_n(\C)$ is identified with
the block matrix
\[
\left(
\begin{array}{cc}
A & \; 0\\
0 & \;1
\end{array}
\right)
\]
in $\GL_{n+1}(\C)$.
Thus $\GL_\infty(\C)$ is the group
of invertible matrices of countable size,
which differ from the identity matrix
only at finitely many entries.

Groups of this form (and related ascending unions
of homogeneous spaces)
have been
considered for a long time in (stable) homotopy theory
and $K$-theory.
Furthermore, results concerning
their representation theory
can be traced back to the 1970s
(more details are given below).
However, only the group structure or topology
was relevant for these studies.
Initially, no attempt was made to consider them
as Lie groups.

The Lie group structure
on $\GL_\infty(\C)$ was first described in~\cite{Mr82}
and that on
$\U_\infty(\C)$ and $\gO_\infty(\R)$
mentioned (cf.\ also page 1053
in the survey article \cite{Mr84}).
The first systematic discussion of direct limits
of finite-dimensional Lie groups was given
in~\cite{NRW1} and \cite{NRW2}. Notably,
a Lie group structure on $G=\dl\,G_n$\vspace{-.4mm}
was constructed there
under technical\linebreak
conditions
which ensure,
in particular, that
\[
\dl\,\exp_{G_n}\colon \dl\,\Lie(G_n)\to \dl\,G_n,\;\;\;
x\mto \exp_{G_n}(x)\;\;
\mbox{for all $\,n\in \N$, $\, x\in \Lie(G_n)$}\vspace{-.6mm}
\]
is a local homeomorphism at~$0$
(see Section~\ref{secconstr}
for the sketch of a more~general construction
from~\cite{FUN}).
Moreover,
situations were described in \cite{NRW2}
where ascending unions of Lie groups (or the corresponding
Lie algebras) can be completed with
respect to some coarser topology.
Also the Lie
group $C^\infty_c(M,H)=\bigcup_{n\in\N}C^\infty_{K_n}(M,H)$
is briefly discussed in \cite{NRW2}
(for finite-dimensional~$H$),
and in~\cite{AH93}.
Test function groups with values
in possibly infinite-dimensional Lie groups
were treated in~\cite{GCX}.
\mbox{Compare} \cite{ACM89} for an early discussion
of gauge groups and automorphism groups
of principal bundles over non-compact
manifolds in the inequivalent
``Convenient Setting''
of infinite-dimensional calculus
(cf.\ also \cite{KaM}).

\noindent
The first construction
of the Lie group structure
on $\Diff_c(M)$ was given in~\cite{Mic},
as part of a discussion
of manifold structures on
spaces of mappings\linebreak
between
non-compact manifolds.
Groups of germs of complex analytic\linebreak
diffeomorphisms of~$\C^n$
around~$0$ were studied in~\cite{Pis}.
The real analytic
analogue was discussed in \cite{Lei},
groups of germs of more general diffeomorphisms
in~\cite{KaR}.
Further recent works will be described later.

It should be stressed that the
current article focusses on
direct limit groups as such, i.e., on their structure and properties.
Representation theory
and harmonic analysis on such groups
are outside its scope.
For completeness, we mention that
the study of (irreducible) unitary representations
of ascending unions of finite
groups started in the 1960s
(see \cite{Th64a} and \cite{Th64b}),
notably for the symmetric group $S_\infty:=\dl\,S_n$.\vspace{-.4mm}
Representations of direct limits
of finite-dimensional
Lie groups were first investigated in the 1970s
(see \cite{Voi} for representations of $\U_\infty(\C)$,
\cite{KS77} for representations of $\U_\infty(\C)$ and $\SO_\infty(\R)$).
Ol'shanski\u{\i} \cite{Ol83}
studied representations for
infinite-dimensional
spherical pairs
like $(\GL_\infty(\R),\SO_\infty(\R))$,
$(\GL_\infty(\C),\U_\infty(\C))$ and
$(\U_\infty(\C),\SO_\infty(\R))$.

The representation theory
of direct limits of both finite groups
and finite-dimensional Lie groups
remains an active area of research.
Representations\linebreak
of
$\gO(\infty,\infty)$, $\U(\infty,\infty)$
and $\Sp(\infty,\infty)$
were studied in \cite{Dv02}
using infinite-dimensional adaptations
of Howe duality.
Novel results concerning the\linebreak
representation theory
of $\U_\infty(\C)$ and $S_\infty$
were obtained in~\cite{Ol03}
and~\cite{KOV}, respectively.
Versions of the Bott-Borel-Weil theorem
for direct limits of finite-dimensional
Lie groups were established in
\cite{NRW4} and \cite{DPW}
(in a more algebraic setting).
J.\,A. Wolf also investigated
principal series\linebreak
representations of suitable
direct limit groups~\cite{Wo05},
as well as the regular\linebreak
representation
on some direct limits of compact
symmetric spaces~\cite{Wo08}.
The paper~\cite{AaK}
discusses
representations of an infinite-dimensional
subgroup of unipotent matrices in $\GL_\infty(\R)$.
Finally, a version of Bochner's theorem for infinite-dimensional
spherical pairs was obtained in~\cite{Ra07}.

There also is a body of literature
devoted to irreducible representations
of diffeomorphism groups of (compact or) non-compact manifolds,
as well as quasi-invariant measures
and harmonic analysis thereon
(see, e.g.,
\cite{VGG75},
\cite{Ki81},
\cite{Hi93},
\cite{Sh01}, and \cite{Sh05}).
Representations of
$C^\infty_c(M,H)$ were studied
by R.\,S. Ismagilov~\cite{Ism}
and in~\cite{AH93}.
Such groups, $\Diff_c(M)$
and semi\-direct products thereof
arise naturally in mathematical physics~\cite{Gol}.

Direct limits of finite-dimensional Lie groups
are also encountered
as dense subgroups of some interesting Banach--Lie groups
(like the group $\U_2(\cH)$
of unitary operators on a complex
Hilbert space~$\cH$ which differ
from $\id_\cH$ by a Hilbert--Schmidt operator)
and other groups of operators.
This frequently enables the calculation
of the homotopy groups of such groups
(see \cite{Pa1}, also \cite{NeH}),
exploiting that the homotopy groups of many direct limits
of classical groups (like $\U_\infty(\C)$)
can be determined using Bott periodicity.
Dense unions of finite-dimensional Lie groups
are also useful in representation theory
(see \cite{Ne98} and \cite{NeO}).

\noindent
We mention a more specialized result:
For very particular classes
of direct limits~$G$ of finite-dimensional
Lie groups, a classification is possible
which uses the
homotopy groups of~$G$ (notably $\pi_1(G)$
and $\pi_3(G)$); see \cite{Ku06}.

In contrast to direct (or inductive) limits,
the dual notion
of an \emph{inverse} (or \emph{projective}) \emph{limit}
of Lie groups was used much earlier
in infinite-dimensional Lie theory.
Omori's theory of ILB-Lie groups
(which are \emph{i}nverse \emph{l}imits
of \emph{B}anach manifolds) gave a strong impetus
to the development of the area in the late 1960s
and early 1970s (see \cite{Om97} and the references therein).
Many important examples of infinite-dimensional
Lie groups could be discussed in this approach,
e.g.\ the group $C^\infty(K,H)=\bigcap_{k\in \N_0}C^k(K,H)$
of smooth maps on a compact manifold~$K$
with values in a finite-dimensional Lie group~$H$,
and the group $\Diff(K)=
\bigcap_{k\in\N}\Diff^k(K)$
of $C^\infty$-diffeomorphisms of a compact\linebreak
manifold.
The passage from compact to non-compact
manifolds naturally leads to the consideration of
direct limits of compactly supported objects.
\section{\,Preliminaries, terminology and basic facts}\label{secprel}
{\bf General conventions.}
We write $\N:=\{1,2,\ldots\}$,
and $\N_0:=\N\cup\{0\}$.
As usual,
$\R$ and $\C$ denote
the fields of real and complex numbers, respectively.
If $(E,\|.\|)$ is a normed space,
$x\in E$ and $r>0$, we write $B^E_r(x):=\{y\in E\colon
\|y-x\|<r\}$.
Topological spaces, topological groups
and locally convex topological vector spaces
are not assumed Hausdorff. However,
manifolds are assumed Hausdorff,
and whenever a locally convex space
serves as the domain or range of a differentiable
map, or as the modelling space
of a Lie group or manifold, it is tacitly
assumed Hausdorff. Moreover, all compact
and all locally compact topological spaces are assumed Hausdorff.
We allow non-Hausdorff topologies
because direct limits are much easier
to describe if the Hausdorff property is omitted
(further explanations will be given at the end of this
section).\\[4mm]
{\bf Infinite-dimensional calculus.}
We are working in the setting of Keller's $C^k_c$-theory~\cite{Kel},
in a topological formulation that avoids the use
of convergence structures
(as in \cite{Mic}, \cite{Mr84}, 
\cite{RES}, \cite{NeS}, and \cite{GaN}).
For more information on
analytic maps, see, e.g., \cite{RES},
\cite{GaN} and (for $\K=\C$)~\cite{BaS}.
\begin{numba}
{\rm Let $\K\in \{\R,\C\}$,
$r\in \N\cup \{\infty\}$,
$E$ and $F$ be locally convex $\K$-vector spaces
and $f\colon U\to F$ be a map
on an open set $U\sub E$.
If $f$ is continuous, we say that $f$ is $C^0$.
We call $f$ a \emph{$C^r_\K$-map}
if $f$ is continuous, the iterated real
(resp., complex) directional
derivatives
\[
d^kf(x,y_1,\ldots, y_k)\; :=\;
(D_{y_k}\cdots D_{y_1}f)(x)
\]
exist for all $k\in \N$ such that $k\leq r$,
$x\in U$ and $y_1,\ldots, y_k\in E$,
and the maps $d^kf\colon U\times E^k\to F$ so obtained
are continuous.
If $\K$ is understood, we write $C^r$ instead
of $C^r_\K$.
If $f$ is $C^\infty$,
we also say that $f$ is \emph{smooth}.
If $\K=\R$,
we say that $f$ is \emph{real analytic}\index{real analytic map}
(or $C^\omega_\R$)
if $f$ extends to a $C^\infty_\C$-map
$\wt{f}\colon \wt{U}\to F_\C$
on an open neighbourhood $\wt{U}$ of~$U$ in the complexification
$E_\C$ of~$E$.}
\end{numba}
\begin{numba}
We mention that a map $f\colon E\supseteq U\to F$ is $C^\infty_\C$
if and only if it is \emph{complex analytic}\index{complex analytic map}
i.e., $f$ is continuous and for each $x\in U$,
there exists a $0$-neighbourhood $Y\sub E$ with $x+Y\sub U$
and continuous homogeneous polynomials $p_n\colon E\to F$
of degree~$n$
such that\vspace{-2mm}
\[
(\forall y\in Y)\;\;\; f(x+y)\, = \sum_{n=0}^\infty \, p_n(y)\,.\vspace{-2mm}
\]
Complex analytic maps are also called \emph{$\C$-analytic}
or $C^\omega_\C$.
\end{numba}
\begin{numba}
It is known that compositions of composable
$C^r_\K$-maps are $C^r_\K$, for each $r\in\N_0\cup\{\infty,\omega\}$.
Thus a
\emph{$C^r_\K$-manifold} $M$ modelled on
a locally convex $\K$-vector space~$E$ can be defined
in the usual way, as a Hausdorff
topological space, together with a maximal
set of homeomorphisms
from open subsets of~$M$ to open subsets of~$E$,
such that the domains cover~$M$ and the transition maps
are $C^r_\K$.
Given $r\in \{\infty,\omega\}$,
a \emph{$C^r_\K$-Lie group} is a group~$G$,
equipped with a structure of $C^r_\K$-manifold
modelled on a locally convex space,
such that the group multiplication and group inversion
are $C^r_\K$-maps. Unless the contrary is stated,
we consider $C^\infty_\K$-Lie groups.
Throughout the following,
the words ``manifold'' and ``Lie group''
will refer to manifolds and Lie groups
modelled on locally convex spaces.
We shall write $T_xM$ for the tangent
space of a manifold $M$ at $x\in M$
and $\Lie(G):=T_{\bf 1}(G)$ for the (topological) Lie algebra
of a Lie group~$G$. Given a $C^1_\K$-map $f\colon M\to N$
between $C^1_\K$-manifolds,
we write $T_xf\colon T_xM\to T_{f(x)}N$ for the tangent map
at $x\in M$. Given a smooth homomorphism $f\colon G\to H$,
we let $\Lie(f):=T_{\bf 1}(f)\colon \Lie(G)\to \Lie(H)$.
\end{numba}
{\bf Direct limits.} \,We recall terminology
and basic facts concerning direct limits.
\begin{numba} (General definitions).
Let $(I,\leq)$ be a directed set,
i.e., $I$ is a non-empty set and $\leq$ a partial order on~$I$
such that any two elements have an upper bound.
Recall that a \emph{direct system}\index{direct system}
(indexed by $(I,\leq)$)
in a category~$\A$ is a pair
$\cS:=((X_i)_{i\in I},
(\phi_{ji})_{j\geq i})$, where each $X_i$ is an object of~$\A$
and $\phi_{ji}\colon X_i\to X_j$
a morphism such that $\phi_{ii}=\id_{X_i}$
and $\phi_{kj}\circ\phi_{ji}=\phi_{ki}$,
for all elements $k\geq j \geq i$ in~$I$.
A \emph{cone over $\cS$}\index{cone}
is a pair $(X,(\phi_i)_{i\in I})$,
where $X$ is an object of~$\A$ and each $\phi_i\colon X_i\to X$
a morphism such that $\phi_j\circ \phi_{ji}=\phi_i$
whenever $j\geq i$. A cone $(X,(\phi_i)_{i\in I})$
is a \emph{direct limit of $\cS$}\index{direct limit}
(and we
write $(X,(\phi_i)_{i\in I})=\dl\,\cS$\vspace{-.3mm}
or $X=\dl\, X_i$),\vspace{-.3mm}
if for every cone $(Y,(\psi_i)_{i\in I})$ over $\cS$,
there
exists a unique morphism $\psi\colon X\to Y$ such that $\psi\circ \phi_i
=\psi_i$
for all $i\in I$.
If ${\cT} = ((Y_i)_{i\in I},(\psi_{ji})_{j\geq i})$ is another
direct system over the same index set, $(Y,(\psi_i)_{i\in I})$
a cone over $\cT$, and $(\eta_i)_{i\in I}$ a family of morphisms
$\eta_i \colon X_i\to Y_i$ which is
\emph{compatible}
with the direct systems
in the sense that $\psi_{ji}\circ \eta_i=\eta_j\circ\phi_{ji}$
for all $j\geq i$, then $(Y,(\psi_i\circ\eta_i)_{i\in I})$
is a cone over $\cS$. We write $\dl\,\eta_i$\vspace{-.9mm}
for the induced morphism $\psi\colon X\to Y$, determined by
$\psi\circ \phi_i=\psi_i\circ \eta_i$.
\end{numba}

Direct limits in the categories of sets,
groups and topological spaces are\linebreak
particularly easy
to understand, and we discuss them now.
Direct limits of\linebreak
topological groups
(which are a more difficult topic)
and direct limits of locally convex spaces
will be discussed afterwards in Sections~\ref{secgp}
and~\ref{seclcx}, respectively.
We concentrate on direct sequences (viz.,
the case $I=\N$)
and actually on ascending sequences, to avoid technical
complications. This is the more justified because
(except for some counterexamples) hardly anything is known about
direct limits of direct systems of Lie groups
which do not admit a cofinal subsequence.
\begin{numba}\label{dlset} (Ascending unions of sets).
If $X_1\sub X_2\sub\cdots$
is an ascending sequence of sets,
let $\phi_{m,n}\colon X_n\to X_m$
be the inclusion map for $m,n\in \N$ with $m\geq n$.
Then $\cS:=((X_n)_{n\in \N},(\phi_{m,n})_{m\geq n}))$
is a direct system in the category $\SET$
of sets and maps.
Define $X:=\bigcup_{n\in \N}X_n$
and let $\phi_n\colon X_n\to X$
be the inclusion map.
Then $(X,(\phi_n)_{n\in \N})$ is a cone over~$\cS$
in $\SET$.
A sequence
of maps $\psi_n\colon X_n\to Y$
to a set~$Y$ gives rise to a cone $(Y,(\psi_n)_{n\in \N})$
if and only if
\[
\psi_m|_{X_n}\;=\; \psi_n\quad\mbox{for all $m,n\in \N$ with $m\geq n$.}
\]
Then $\psi\colon X\to Y$, $\psi(x):=\psi_n(x)$ if $n\in \N$
and $x\in X_n$ is a well-defined map,
and is uniquely determined by the requirement
that $\psi\circ \phi_n=\psi|_{X_n}=\psi_n$ for each $n\in\N$.
Thus $(X,(\phi_n)_{n\in \N})=\dl\,\cS$ in $\SET$.
\end{numba}
\begin{numba}\label{pardlgp}
(Direct limits of groups).
If each $X_n$ is a group in the situation of \ref{dlset}
and each $\phi_{m,n}$ a homomorphism,
then $\cS$ is a direct system in the category~$\G$
of groups and homomorphisms.
If $x,y\in X$, there exists $n\in \N$
such that $x,y\in X_n$. We define the product of
$x$ and~$y$ in~$X$ as their product in~$X_n$,
i.e., $x\cdot y =\phi_n(x)\cdot \phi_n(y):=\phi_n(x\cdot y)$.
Since each $\phi_{m,n}$ is a homomorphism, $x\cdot y$
is independent of the choice of~$n$,
and it is clear that the product so defined makes~$X$
a group and each $\phi_n\colon X_n\to X$ a homomorphism.
If $(Y,(\psi_n)_{n\in\N})$ is a cone over~$\cS$ in~$\G$,
let $\psi\colon X\to Y$ be the unique map
such that $\psi\circ \phi_n=\psi_n$ for each $n\in \N$,
as in~\ref{dlset}.
Given $x,y\in X$, say $x,y\in X_n$, we then have
$\psi(xy)=\psi(\phi_n(xy))=\psi_n(xy)=\psi_n(x)\psi_n(y)=\psi(x)\psi(y)$,
whence~$\psi$ is a\linebreak
homomorphism.
Thus $(X,(\phi_n)_{n\in \N})=\dl\,\cS$
in~$\G$.\\[1mm]
If $\cS=((X_n)_{n\in \N},(\phi_{m,n}))$
is a direct sequence of groups
(with $\phi_{m,n}\colon X_n\to X_m$
not necessarily injective), then $K_n:=\bigcup_{m\geq n}\ker(\phi_{m,n})$
is a normal subgroup of $X_n$. Consider the quotient groups
$G_n:=X_n/K_n$, the canonical quotient maps
$q_n\colon X_n\to G_n$
and the homomorphisms $\psi_{m,n}\colon G_n\to G_m$
determined by $\psi_{m,n}\circ q_n=q_m\circ \phi_{m,n}$.
Then each $\psi_{m,n}$ is injective
and it is clear that the direct limit $(G,(\psi_n)_{n\in \N})$
of the ``injective quotient system''
$((G_n)_{n\in\N},(\psi_{m,n}))$
yields a direct limit $(G,(\psi_n\circ q_n)_{n\in \N})$
of~$\cS$ (cf.\ \cite[\S3]{NRW1}).
\end{numba}
\begin{numba}\label{DLtop}
(Direct limits of topological spaces).
If each $X_n$ is a topological space
in the situation of~\ref{dlset}
and each $\phi_{m,n}\colon X_n\to X_m$
a continuous map, we equip $X=\bigcup_{n\in \N}X_n$
with the finest topology $\cO_{\DL}$ making each inclusion map
$\phi_n\colon X_n\to X$ continuous
(the so-called
\emph{direct limit topology}).\index{direct limit topology}
Thus $U\sub X$ is open (resp., closed)
if and only if $\phi_n^{-1}(U)=U\cap X_n$ is open (resp., closed)
in~$X_n$ for each $n\in \N$.
Then $(X,(\phi_n)_{n\in \N})=\dl\,\cS$\vspace{-.7mm} in
the category $\TOP$ of topological spaces and continuous maps.
To see this, let
$(Y,(\psi_n)_{n\in \N})$ be a cone over~$\cS$ in~$\TOP$.
Let $\psi\colon X\to Y$ be the unique map with
$\psi\circ \phi_n=\psi_n$ for each~$n$.
If $U\sub Y$ is open, then
$\psi^{-1}(U)\cap X_n=(\psi|_{X_n})^{-1}(U)=(\psi_n)^{-1}(U)$
is open in~$X_n$, for each~$n$.
Hence $\psi^{-1}(U)$ is open in~$X$
and thus $\psi$ is continuous.\\[2.5mm]
The direct system~$\cS$ is called
\emph{strict}\index{strict direct system}
if each $\phi_n$ is a topological embedding
(i.e., $X_{n+1}$ induces the topology of~$X_n$).
Then also the inclusion map $\phi_n\colon X_n\to X$ is a
topological embedding for each~$n$ \cite[Lemma~A.5]{NRW2}.
It is also known that~$X$ has the separation property $T_1$
if each $X_n$ is $T_1$ (see, e.g., \cite[Lemma~1.7\,(a)]{FUN}).
And in the case of a direct sequence
$X_1\sub X_2\sub\cdots$ of\linebreak
locally compact
spaces~$X_n$, the direct limit topology on $\bigcup_{n\in \N}X_n$
is Hausdorff (as observed in
\cite[Lemma~1.7\,(c)]{FUN},
the strictness hypotheses
in \cite[Proposition~4.1\,(ii)]{Han}
and \cite[Lemma~3.1]{DIR}
is unnecessary).
\end{numba}
{\bf First remarks on ascending unions of Lie groups and direct limits.}
Consider an ascending sequence $G_1\sub G_2\sub\cdots$
of $C^\infty_\K$-Lie groups,
such that the inclusion maps $j_{m,n}\colon G_n\to G_m$
are $C^\infty$-homomorphisms for all $m,n\in \N$ with $m\geq n$.
Then $\cS:=((G_n)_{n\in \N},(j_{m,n})_{m\geq n})$
is a direct system in the category $\LIE_\K$
of $C^\infty_\K$-Lie groups
and $C^\infty_\K$-homomorphisms.
One would not expect
that $\cS$ always
has a direct limit in the category
of $C^\infty_\K$-Lie groups
(although no counterexamples
are known at the time of writing).
What is more,
there is no general
construction principle
for a Lie group structure on $\bigcup_{n\in \N}G_n$
such that all inclusion maps $j_n\colon G_n\to G$
are $C^\infty_\K$-homomorphisms
(unless restrictive
conditions are imposed,
as in Section~\ref{secconstr}).
\begin{numba}\label{explainnow}
However,
in many concrete cases
we are given
such a Lie group structure on $G:=\bigcup_{n\in \N}G_n$.
Then $(G,(j_n)_{n\in \N})$ is a cone over~$\cS$
in $\LIE_\K$,
and it is natural to ask if $G=\dl\, G_n$\vspace{-.7mm}
as a Lie group. A sequence
of $C^\infty_\K$-homomorphisms $f_n\colon G_n\to H$
to a $C^\infty_\K$-Lie group~$H$
is a cone over $\cS$ if and only if
\[
f_m|_{G_n}\;=\; f_n\quad\mbox{for all $m,n\in \N$ with $m\geq n$.}
\]
Then $f\colon G\to H$, $f(x):=f_n(x)$ if $n\in \N$
and $x\in G_n$ is a well-defined
homo\-mor\-phism. This map
is uniquely determined by the requirement
that $f\circ j_n=f|_{G_n}=f_n$ for each $n\in\N$.
Therefore, $(G,(j_n)_{n\in \N})=\dl\,\cS$
holds in $\LIE_\K$
if and only if each $f$ of the preceding form
is~$C^\infty_\K$. A similar argument applies if~$H$
is a topological group, smooth manifold or topological space.\\[2.5mm]
Thus questions (a)--(d)
posed in \ref{dlpropannounce}
amount to asking if
$G=\dl\, G_n$\vspace{-.5mm} holds\,\ldots\vspace{-1mm}
\begin{description}[(a)$'$]
\item[(a)$'$]
in the category of smooth manifolds
(modelled on locally convex spaces)
and smooth maps between them?
\item[(b)$'$]
in the category of topological spaces
and continuous maps?
\item[(c)$'$]
in the category $\LIE_\K$ of Lie groups?
\item[(d)$'$]
in the category of topological groups
and continuous homomorphisms?
\end{description}
\end{numba}
{\bf The Hausdorff property.}
We allow non-Hausdorff topologies
because direct limits are much easier to
describe if the Hausdorff property is omitted.
In fact, we have already seen
that it is always possible to
topologize a union $X=\bigcup_{n\in \N}X_n$
of topological spaces
in such a way that it becomes
the direct limit $\dl\,X_n$\vspace{-.7mm}
in the category of topological spaces
(see \ref{DLtop}),
and likewise
a union of topological groups
(resp., locally convex spaces)
can always be made the direct limit
in the category of topological groups
resp., locally convex spaces
(see Sections~\ref{secgp}
and~\ref{seclcx}). A mere union $X=\bigcup_{n\in \N}X_n$
is a very concrete object, and easy to work with.

By contrast, if each $X_n$ is Hausdorff,
then the direct limit $\dl\,X_n$\vspace{-.7mm}
in the category of Hausdorff
topological spaces (resp., Hausdorff topological
groups, resp., Hausdorff locally convex spaces)
can only be realized as a
quotient of $X=\bigcup_{n\in \N}X_n$ in general,
and is a much more elusive object in this case.

Luckily, in all situations we are interested in,
$X$ from above injects\linebreak
continuously into a Lie group
and thus $X$ is Hausdorff. Then automatically~$X$
also is the direct limit in the category
of Hausdorff topological spaces (resp., Hausdorff topological
groups, resp., Hausdorff locally convex spaces).
\section{\,Direct limits of topological groups}\label{secgp}
As an intermediate step towards the study of Lie groups,
let us consider a sequence
$G_1\sub G_2\sub\cdots$ of topological groups,
such that all inclusion maps $G_n\to G_{n+1}$
are continuous homomorphisms.
We make $G=\bigcup_{n\in \N}G_n$
the direct limit group (as in \ref{pardlgp})
and give it the finest group
topology $\cO_{\DLG}$ making each inclusion map
$G_n\to G$ continuous. Then $G=\dl\,G_n$\vspace{-.7mm}
in the category of (not necessarily Hausdorff)
topological groups. Moreover,
if each $G_n$ is Hausdorff,
then the factor group of~$G$ modulo
the closure $\wb{\{{\bf 1}\}}\sub G$ is the direct limit
in the category of Hausdorff topological groups.\footnote{If $G$ is Hausdorff,
then no passage to the quotient
is necessary.}

Unfortunately, the preceding description of the
topology $\cO_{\DLG}$ on the\linebreak
direct limit topological group is not at all concrete.
Various questions are natural (and also
relevant for our studies of Lie groups):
Does $\cO_{\DLG}$ coincide with the direct
limit topology $\cO_{\DL}$ (as in \ref{DLtop})?
\,Can $\cO_{\DLG}$ be described more explicitly?
Given a group topology on $G=\bigcup_{n\in \N}G_n$,
how can we prove that it agrees with $\cO_{\DLG}$?

We now give some answers to the first and last
question. An answer to the second question,
namely the description of $\cO_{\DLG}$
as a so-called ``bamboo-shoot'' topology,
can be found in~\cite{TSH} and \cite{HSTH}
(under suitable\linebreak
hypotheses).\\[4mm]
{\bf Comparison of {\boldmath$\cO_{\DL}$} and {\boldmath$\cO_{\DLG}$}.}
It is clear from the definition that the\linebreak
direct limit topology $\cO_{\DL}$ is finer
than $\cO_{\DLG}$.
Moreover, $\cO_{\DL}$ may be properly finer
than $\cO_{\DLG}$,
as emphasized by Tatsuuma et al.\ \cite{TSH}.\footnote{In part of the
older literature, there was some confusion concerning
this point.}
To understand this difficulty, let $\eta_n\colon G_n\to G_n$, $x\mto x^{-1}$
and $\eta\colon G\to G$ be the inversion maps
and $\mu_n\colon G_n\times G_n\to G_n$, $(x,y)\mto xy$
as well as $\mu\colon G\times G\to G$
be the respective group multiplication.
Then
\[
\eta\;=\; \dl\,\eta_n\colon \big(\dl\, G_n, \cO_{\DL}\big)\to \big(\dl\, G_n,\cO_{\DL}\big)
\]
is always continuous.
However, it may happen that $\mu$ is discontinuous
(with respect to the product topology on $G\times G$),
in which case $(G,\cO_{\DL})$ is not a topological
group and hence $\cO_{\DL}\not=\cO_{\DLG}$.
We recall a simple example for this pathology
from~\cite{TSH}:
\begin{example}\label{badDL1}
Let $G_n:=\Q\times \R^{n-1}$
with the addition and topology induced\linebreak
by $\R^n$.
Identifying $\R^{n-1}$ with the vector subspace $\R^{n-1}\! \times \{0\}$
of $\R^n$, we
obtain a strict direct sequence
$G_1\sub G_2\sub\cdots$
of metrizable topological groups.
It can be shown by direct calculation
that the direct limit topology
$\cO_{\DL}$
does not make the group multiplication
on $G:=\bigcup_{n\in \N}G_n$ continuous (see \cite[Example~1.2]{TSH}).
\end{example}
To understand the difficulties concerning the
group multiplication (in\linebreak
contrast to the group inversion)
on $G=\bigcup_{n\in \N}G_n$, note that we always
have a continuous map
\[
\dl\,\mu_n\colon \big(\dl\, (G_n\times G_n), \cO_{\DL}\big)\to \big(\dl\, G_n,\cO_{\DL}\big)\,.
\]
Thus $\mu$ is continuous as a map
from $(G\times G, \cO_{\DL})$ to $(G,\cO_{\DL})$,
i.e., it becomes continuous if,
instead of the product
topology, the topology $\cO_{\DL}$ is used on $G\times G$
which makes it the direct limit topological space
$\dl\,(G_n\times G_n)$.\vspace{-.7mm}
This topology is finer than the product topology
and, in general, properly finer.
If the direct limit topology
on $G\times G$ happens to coincide
with the product topology,
then $(G,\cO_{\DL})$
is a topological group and thus $\cO_{\DL}=\cO_{\DLG}$
(cf.\ \cite{HSTH} and \cite[\S3]{DIR}).
The following proposition describes a
situation where the two topologies
coincide. We
recall that a topological
space $X$ is said to be a
\emph{$k_\omega$-space}\index{$k_\omega$-space}
if it is the direct limit topological space
of an ascending sequence
$K_1\sub K_2\sub\cdots$ of compact topological
spaces (see, e.g., \cite{GGH} and the references
therein).\footnote{These spaces can also be characterized as
the hemicompact $k$-spaces.}
Such spaces are always Hausdorff (see \ref{DLtop}).
For example, every $\sigma$-compact, locally compact
space is a $k_\omega$-space.
A topological space~$X$ is
called \emph{locally~$k_\omega$}\index{locally $k_\omega$ space}
if every point $x\in X$ has an open neighbourhood in~$X$
which is a $k_\omega$-space
in the induced topology \cite[Definition~4.1]{GGH}.
E.g., every locally compact topological
space is locally~$k_\omega$.
The topological space underlying
a topological group~$G$ is locally~$k_\omega$
if and only if~$G$ has an open subgroup
which is a $k_\omega$-space \cite[Proposition~5.3]{GGH}.
See \cite[Proposition~4.7]{GGH}
for the following fact.
The special case where each $X_n$ and $Y_n$
is locally compact was first proved in \cite[Theorem~4.1]{HSTH}
(cf.\ also \cite[Proposition~3.3]{DIR} for the strict case):
\begin{proposition}
Let
$X_1\sub X_2\sub \cdots$
and $Y_1\sub X_2\sub \cdots$
be topological spaces
with continuous inclusion maps $X_n\to X_{n+1}$
and $Y_n\to Y_{n+1}$.
If each $X_n$ and each $Y_n$
is locally $k_\omega$,
then
\[
\dl\, (X_n\times Y_n)\;=\; \big(\dl\,X_n\big)\times\big(\dl\,Y_n\big)
\]
as a topological space.
\end{proposition}
Using that direct limits of ascending sequences
of locally $k_\omega$-spaces are\linebreak
locally~$k_\omega$
by \cite[Proposition~4.5]{GGH}
(and thus Hausdorff),
the preceding discussion
immediately entails the following conclusion 
from~\cite{GGH}
(cf.\ \cite[Theorem~2.7]{TSH}
for locally compact $G_n$,
as well as \cite[Corollary~3.4]{DIR}
(in the case of a strict direct system)).
\begin{corollary}
Consider a sequence
$G_1\sub G_2\sub\cdots$
of topological groups
such that each inclusion map $G_n\to G_{n+1}$
is a continuous homomorphism.
If the topological space underlying $G_n$
is locally $k_\omega$ for each $n\in \N$
$($for example, if each $G_n$ is locally compact$)$,
then the direct limit topology is Hausdorff
and makes $G=\bigcup_{n\in \N}G_n$
the direct limit topological group.
\end{corollary}
\begin{numba}
Given topological groups
$G_1\sub G_2\sub \cdots$
such that all inclusion maps $G_n\to G_{n+1}$
are continuous homomorphisms,
consider the conditions:
\begin{description}[(D)]
\item[(a)]
$G_n$ is an open subgroup of $G_{n+1}$
(with the induced topology)
for all sufficiently large $n$.
\item[(b)]
For each sufficiently large~$n$,
the topological group $G_n$ has an identity neighbourhood $U$
whose closure in $G_m$ is compact
for some $m\geq n$.
\end{description}
Then $\cO_{\DL}=\cO_{\DLG}$ holds
if (a) or (b) is satisfied~\cite[Theorems~2 and 3]{Yam}.
By a most remarkable theorem of Yamasaki\index{Yamasaki's Theorem}
\cite[Theorem~4]{Yam},
the validity of (a) or (b)
is also \emph{necessary} in order that $\cO_{\DL}=\cO_{\DLG}$,
provided that each $G_n$ is metrizable
and the inclusion maps $G_n\to G_{n+1}$
are topological embeddings.
\end{numba}
{\bf Criteria ensuring that a given group topology
coincides with {\boldmath$\cO_{\DLG}$}.}
Frequently, a given
topological group~$G$
is a union $G=\bigcup_{n\in \N}G_n$
of topological groups,
such that all inclusion maps $G_n\to G_{n+1}$
and $G_n\to G$ are continuous homomorphisms.
In many cases,
a criterion from~\cite{COM}
helps to see that
the given topology on~$G$ coincides
with $\cO_{\DLG}$ (cf.\ \cite[Proposition~11.8]{COM}).\\[2.5mm]
The criterion uses the weak direct product
$\prod_{n\in \N}^*G_n$\index{weak direct product} as a tool.
The latter can be formed for
any sequence $(G_n)_{n\in \N}$ of topological
groups. It is defined as the subgroup
of all $(g_n)_{n\in \N}\in \prod_{n\in \N}G_n$
such that $g_n=1$ for all but finitely many~$n$.
The weak direct product is a topological
group; a basis for its topology (the so-called
``box topology'')\index{box topology}
is given by sets of the form
$\prod_{n\in \N}U_n\cap \prod^*_{n\in\N}G_n$
(the ``boxes''),\index{box}
where $U_n\sub G_n$ is open for each $n$
and ${\bf 1}\in U_n$ for almost all~$n$.\\[2.5mm]
Returning to the case where $G_1\sub G_2\sub\cdots$
and $G=\bigcup_{n\in \N}G_n$, we can consider the ``product
map''\index{product map}
\[
\pi\colon
{\textstyle \prod_{n\in \N}^*} \, G_n \to G\,,
\quad (g_n)_{n\in \N} \mto g_1g_2\cdots g_N\,,
\]
where $N\in \N$ is so large that $g_n={\bf 1}$ for all $n>N$.
\begin{proposition}\label{propcriter}
If the product map $\pi\colon
\prod_{n\in \N}^* G_n \to G$ is open at~${\bf 1}$,
then the given topology on~$G$ coincides with $\cO_{\DLG}$
and thus $G=\dl\, G_n$\vspace{-.3mm}
as a topological group.
The openness of $\pi$ at~${\bf 1}$ is guaranteed
if there exists a map $\sigma\colon \Omega \to \prod^*_{n\in \N}G_n$
on an identity neighbourhood $\Omega\sub G$
such that $\pi\circ\sigma=\id_\Omega$,
$\sigma({\bf 1})={\bf 1}$ and $\sigma$ is continuous at~${\bf 1}$.
\end{proposition}
\begin{remark}
Such a section~$\sigma$ to $\pi$ might be called
a \emph{fragmentation map},\index{fragmentation map}
in analogy to concepts in the theory
of diffeomorphism groups (cf.\ \cite[\S2.1]{Ban}).
\end{remark}
\begin{example}
It can be shown that the Lie groups
$\Diff_c(M)=\bigcup_{n\in \N}\Diff_{K_n}(M)$
and $C^r_c(M,H)=\bigcup_{n\in \N}C^r_{K_n}(M,H)$
(as defined in the introduction)
always admit fragmentation maps
(even smooth ones);
cf.\ \cite[Lemmas~5.5 and 7.7]{COM}.
Hence
$\Diff_c(M)=\dl\,\Diff_{K_n}(M)$
and $C^r_c(M,H)=\dl\,C^r_{K_n}(M,H)$\vspace{-.7mm}
as topological groups.
\end{example}
\section{\,Non-linear mappings on locally convex
direct limits}\label{seclcx}
Consider a sequence $E_1\sub E_2\sub \cdots$
of locally convex spaces,
such that each inclusion map $E_n\to E_{n+1}$
is continuous and linear.
Then there is a finest locally convex vector topology $\cO_{\lcx}$
on $E:=\bigcup_{n\in \N}E_n$
making each inclusion map $E_n\to E$
continuous,
called the \emph{locally convex direct limit
topology}.\footnote{This
topology can be described in various ways.
We mention: (1) A convex set $U\sub E$ is
open if and only if $U\cap E_n$ is open in~$E_n$,
for each~$n\in \N$.
(2) A seminorm $q\colon E\to [0,\infty[$ is
continuous if and only if $q|_{E_n}$ is
continuous, for each $n\in \N$.}\index{locally convex direct limit topology}
\begin{numba}\label{baslcx}
Some basic properties
of locally convex direct limits
are frequently used:
\begin{description}[(D)]
\item[(a)]
If the direct sequence $E_1\sub E_2\sub \cdots$
is strict, then
$(E,\cO_{\lcx})$ induces the given topology on~$E_n$,
for each $n\in \N$
(see Proposition~9\,(i) in
\cite[Chapter~II, \S4, no.\,6]{BTV}).
\item[(b)]
If the direct sequence $E_1\sub E_2\sub \cdots$
is strict and each $E_n$ is Hausdorff,
then also $(E,\cO_{\lcx})$ is Hausdorff
(see Proposition~9\,(i) in
\cite[Chapter~II, \S4, no.\,6]{BTV}).
\item[(c)]
If the direct sequence $E_1 \sub E_2\sub \cdots$
is strict and each $E_n$ complete,
then the locally convex direct limit $E=\bigcup_{n\in \N}E_n$
is boundedly regular
(cf.\ Proposition~6
in \cite[Chapter~III, \S1, no.\,4]{BTV})
and hence also compactly regular, in view of~(a).
\item[(d)]
If the direct sequence $E_1\sub E_2\sub \cdots$ is strict and each $E_n$
complete, then also the locally
convex direct limit $E$ is complete
(see Proposition~9\,(iii) in
\cite[Chapter~II, \S4, no.\,6]{BTV}).
\item[(e)]
If also $F_1\sub F_2\sub\cdots$ is an ascending sequence
of locally convex spaces, with locally convex
direct limit
$F=\bigcup_{n\in \N}F_n$,
then the locally convex\linebreak
direct
limit topology on $\bigcup_{n\in \N}(E_n\times F_n)$
and the product topology on $E\times F$ coincide
\cite[Theorem~3.4]{HSTH}
(because finite direct products
coincide with finite direct sums
in the category of locally convex spaces).
\end{description}
\end{numba}

The reader may find \cite{Fl80}
and \cite{Bie} convenient points of entry to
the research literature on locally convex direct limits.

We mention that
few general results ensuring the
Hausdorff property for
locally convex direct limits
$E=\dl\, E_n$\vspace{-.7mm}
are known
(besides \ref{baslcx}\,(b) just encountered
and Proposition~\ref{propslv} below).
Some Hausdorff criteria for
direct limits of
Banach spaces (and normed spaces) can be found in \cite{Fl79}
(see also \cite[p.\,214]{Fl80}).
In concrete examples, a very simple argument
frequently works:
If one can find an injective continuous linear map
from $E$ to a some Hausdorff locally convex space,
then~$E$ is Hausdorff.
However, the fact remains that
non-Hausdorff locally convex
direct limits do exist:
See~\cite{Mk63} for examples
where each $E_n$ is a Banach space;
\cite[p.\,207]{Fl80} for a simple example (due to \mbox{L.~Waelbroeck)}
with each $E_n$ a normed space;
and \cite[p.\,227, Corollary~2]{Fl80} for an example where
each $E_n$ is a nuclear Fr\'{e}chet space (cf.\ also \cite{Sm69}).

It is well known that $\cO_{\lcx}$ and $\cO_{\DLG}$
coincide on $E=\bigcup_{n\in \N}E_n$,
because on $\bigoplus_{n\in \N}E_n=\bigcup_{n\in \N}
\prod_{k=1}^nE_k$,
both $\cO_{\lcx}$ and $\cO_{\DLG}$
coincide with the box topology
and $\bigcup_{n\in \N}E_n$
(with either topology) can be considered
as a quotient of the direct sum
(see \cite[Lemma~2.7]{COM}; cf.\ \cite[Proposition~3.1]{HSTH}
for a different argument. Cf.\ also
\cite{Ko69} and \cite[Chapter~II, Exercise 14 to \S4]{BTV}).

It is also known that $\cO_{\lcx}$ need
not coincide with $\cO_{\DL}$ (see \cite{Sr59}
or\linebreak
Exercise 16\,(a) to \S4 in \cite[Chapter~II]{BTV};
cf.\ also \cite[p.\,506]{Dud}).
E.g., $\cO_{\lcx}$ is properly coarser than
$\cO_{\DL}$ if each $E_n$ is an infinite-dimensional
Fr\'{e}chet space and $E_n$ is a proper vector subspace of $E_{n+1}$
with the induced topology, for each $n\in \N$
(\cite[Proposition~4.26\,(ii)]{KaM};
cf.\ Yamasaki's Theorem recalled in
Section~\ref{secgp}).
The following concrete
example shows that not even
smoothness or analyticity
of $f|_{E_n}$
ensures that a map
$f\colon E\to F$
on a locally convex direct limit $E=\bigcup_{n\in \N}E_n$
is continuous (let alone smooth or analytic).
\begin{example}\label{extenso}
Consider the map
\[
g\colon C^\infty_c(\R,\C)\to C^\infty_c(\R\times \R,\C)\,,\quad
g(\gamma)\,:=\, \gamma\otimes \gamma
\]
between spaces of compactly supported
smooth functions,
where $(\gamma\otimes \gamma)(x,y)$
$:=\gamma(x)\gamma(y)$
for $x,y\in \R$. It can be shown that
$g$ is discontinuous, although
$g|_{C^\infty_{[{-n},n]}(\R,\C)}\colon
C^\infty_{[{-n},n]}(\R,\C)\to C^\infty_c(\R\times \R,\C)$
is a continuous homogeneous polynomial
(and hence complex analytic), for each $n\in \N$
(see Remark~7.9 in \cite{COM},
based on \cite[Theorem~2.4]{HSTH}).
\end{example}
\begin{remark}\label{whybad}
Consider the locally convex direct limit
$E=\bigcup_{n\in \N}E_n$
of Hausdorff locally convex spaces $E_1\sub E_2\sub \cdots$
over $\K\in \{\R,\C\}$.
Let
$U_1\sub U_2\sub \cdots$ be an ascending sequence
of open sets $U_n\sub E_n$,
and $U:=\bigcup_{n\in \N}U_n$.
Let $r\in \N\cup \{\infty\}$,
$F$ be a Hausdorff locally convex space
and $f\colon U\to F$ be a map such that
$f|_{U_n}\colon E_n\supseteq U_n\to F$ is $C^r_\K$
for each $n\in \N$.
Assume that $E$ is Hausdorff and $U\sub E$ is open.\footnote{E.g.,
we might start with an open set $U\sub E$ and set
$U_n:=U\cap E_n$.}
Then the iterated directional
derivatives
\[
d^kf(x,y_1,\ldots,y_k)\;=\;
(D_{y_k}\cdots D_{y_1}f)(x)
\]
exist for all $k\in \N$ with $k\leq r$
and all $x\in U$ and $y_1,\ldots, y_k\in E$,
because $x\in U_n$ and $y_1,\ldots, y_k\in E_n$
for some $n\in \N$ and then
$(D_{y_k}\cdots D_{y_1}f)(x)=d^k(f|_{U_n})(x,y_1,\ldots, y_k)$.
Hence only continuity
of the maps $d^kf$,
which satisfy
\begin{equation}\label{givesCk}
d^kf|_{U_n\times (E_n)^k}\;=\; d^k(f|_{U_n})\quad\mbox{for all $n\in \N$,}
\end{equation}
may be missing for some~$k$,
and may prevent $f$ from being a $C^r_\K$-map.
\end{remark}
We mention that locally convex direct limits
of ascending sequences of Banach spaces
(resp., Fr\'{e}chet spaces) are called
(LB)-spaces (resp., (LF)-spaces).
If the sequence is strict,
we speak of LB-spaces (resp., LF-spaces).\footnote{These\hspace*{-.1mm}
conventions\hspace*{-.1mm} are\hspace*{-.1mm} local.\hspace*{-.1mm}
The\hspace*{-.1mm} meanings\hspace*{-.1mm} of\hspace*{-.1mm}
`LF'\hspace*{-.1mm} and\hspace*{-.1mm}
`(LF)'\hspace*{-.1mm}
vary\hspace*{-.1mm} in\hspace*{-.1mm}
the\hspace*{-.1mm} literature.}\index{(LB)-space}\index{LB-space}
\index{(LF)-space}\index{LF-space}
A locally convex space~$E$ is called
a \emph{Silva space}\index{Silva space}
if it is the locally convex direct limit
of an ascending sequence $E_1\sub E_2\sub\cdots$
of Banach spaces, such that all inclusion maps
$E_n\to E_{n+1}$ are compact operators
(cf.\ \cite{Si55} and \cite{Fl71}).\footnote{A locally
convex space is a Silva space
if and only if it is isomorphic to the dual
of a Fr\'{e}chet-Schwartz space~\cite{Fl71};
therefore Silva spaces are also called
(DFS)-spaces.}

Silva spaces are
very well-behaved
direct limits. We recall from~\cite{Fl71}:
\begin{proposition}\label{propslv}
If $E=\bigcup_{n\in \N}E_n$ is a Silva space, then
the following hold:
\begin{description}[(D)]
\item[\rm(a)]
$E$ is Hausdorff and complete;
\item[\rm(b)]
$E=\bigcup_{n\in \N}E_n$ is boundedly regular
and hence also compactly regular;\,\footnote{Using that
the inclusion maps $E_n\to E_{n+1}$ are
compact operators.}
\item[\rm(c)]
The locally convex direct limit topology on~$E$
coincides with the direct limit topology $\cO_{\DL}$;
\item[\rm(d)]
If also $F=\bigcup_{n\in \N}F_n$ is a Silva space,
with $F_n\to F_{n+1}$ compact,
then $E\times F=\bigcup_{n\in \N}(E_n\times F_n)$
is a Silva space.\footnote{The inclusions
$E_n\times F_n\to E_{n+1}\times F_{n+1}$ are
compact operators, and \ref{baslcx}\,(e) holds.}
\end{description}
\end{proposition}
Some interesting
infinite-dimensional Lie groups are modelled
on Silva spaces, e.g.\ the group
$\Diff^\omega(K)$ of real analytic diffeomorphisms
of a compact real analytic manifold~$K$ (see \cite{Les};
cf.\ \cite[Theorem~43.4]{KaM}).
More examples will be encountered below.\\[3mm]
{\bf Mappings on Silva spaces or unions of
{\boldmath $k_\omega$}-spaces.}
In good cases,
the pathology described in Remark~\ref{whybad}
cannot occur (see \cite[Lemma~9.7]{COM}
and \cite[Proposition~8.12]{GGH}):
\begin{proposition}\label{silvkom}
Consider the locally convex direct limit
$E=\bigcup_{n\in \N}E_n$
of Hausdorff locally convex spaces $E_1\sub E_2\sub \cdots$
over $\K\in \{\R,\C\}$.
Let\linebreak
$U_1\sub U_2\sub \cdots$ be an ascending sequence
of open sets $U_n\sub E_n$,
and $U:=\bigcup_{n\in \N}U_n$.
Let $r\in \N_0\cup \{\infty\}$,
$F$ be a Hausdorff locally convex space
and $f\colon U\to F$ be a map such that
$f|_{U_n}$ is $C^r_\K$
for each $n\in \N$.
Assume that
\begin{description}[(D)]
\item[\rm(a)]
Each $E_n$ is a $k_\omega$-space; or:
\item[\rm(b)]
$E_n$ is a Banach space
and the inclusion map $E_n\to E_{n+1}$
a compact operator, for each $n\in \N$
$($in which case $E$ a Silva space$)$.
\end{description}
Then $E$ is Hausdorff
and the locally convex direct
limit topology on~$E$ coincides with $\cO_{\DL}$.
Moreover, $U$ is open in~$E$
and $f\colon U\to F$ is $C^r_\K$.
\end{proposition}
In the Silva case, the hypotheses
of Proposition~\ref{silvkom}
can be relaxed
(cf.\ \cite[Proposition~2.8]{Ls85}).
Real analyticity is more elusive.
E.g., there exists a real-valued
map~$f$ on the Silva space $\R^{(\N)}:=\dl\,\R^n$\vspace{-.5mm}
which is not real analytic
although $f|_{\R^n}$ is real analytic
for each $n\in \N$ (cf.\ \cite[Example~10.8]{KaM}).\\[3mm]
{\bf Complex analytic maps on (LB)-spaces.}
A very useful result
from \cite{Dah}
frequently facilitates
to check complex analyticity
beyond Silva spaces.
\begin{theorem}[Dahmen's Theorem]\label{thmdah}
Let $E_1\sub E_2\sub\cdots$\index{Dahmen's Theorem}
be an ascending sequence of normed spaces
$(E_n,\|.\|_n)$ over~$\C$ such that, for each $n\in \N$,
the inclusion map
$E_n\to E_{n+1}$ is continuous and complex linear,
of operator norm at most~$1$.
Let $r\in \;]0,\infty[$,
$U_n:=\{x\in E_n\colon \|x\|_n<r\}$
for $n\in \N$,
and $F$ be a complex
locally convex space.
Assume that the locally convex direct
limit $E=\bigcup_{n\in \N}E_n$ is Hausdorff.
Then $U:=\bigcup_{n\in \N}U_n$
is open in~$E$ and if
\[
f\colon U \to F
\]
is a map
such that $f|_{U_n}\colon E_n\supseteq U_n\to F$
is complex analytic and bounded
for each $n\in \N$,
then $f$ is complex analytic.
\end{theorem}
{\bf Mappings between direct sums.} If $(E_n)_{n\in \N}$
is a sequence of locally\linebreak
convex spaces,
we equip
$\bigoplus_{n\in \N} E_n$
with the box topology (as introduced before
Proposition~\ref{propcriter}).
See \cite[Proposition~7.1]{MEA} for the following result.
\begin{proposition}\label{dsummp}
Let $(E_n)_{n\in \N}$ and $(F_n)_{n\in \N}$
be sequences of Hausdorff\linebreak
locally convex spaces,
$r\in \N_0\cup\{\infty\}$,
$U_n\sub E_n$ be open
and $f_n\colon U_n\to F_n$ be~$C^r$.
Assume that $0\in U_n$
and $f_n(0)=0$ for all but finitely many
$n\in \N$.
Then $\bigoplus_{n\in \N}U_n:=(\bigoplus_{n\in \N} E_n)\cap \prod_{n\in\N}U_n$
is open in $\bigoplus_{n\in\N}E_n$ and the
map
$\oplus_{n\in \N}f_n\colon \bigoplus_{n\in \N}U_n\to
\bigoplus_{n\in \N}F_n$,
$(x_n)_{n\in \N}\mto (f_n(x_n))_{n\in \N}$
is $C^r$.\vspace{1mm}
\end{proposition}
{\bf Non-linear maps between spaces of test functions.}
Let \mbox{$r,s \in \N_0\cup\{\infty\}$,}
$M$ be a $\sigma$-compact, finite-dimensional $C^r$-manifold,
$N$ be a $\sigma$-compact,
finite-dimensional $C^s$-manifold,
$E$, $F$ be Hausdorff locally convex spaces,
$\Omega \sub C^r_c(M,E)$ be open
and $f\colon \Omega\to C^s_c(N,F)$ be a map.
We say that $f$ is \emph{almost local}\index{almost local map}
if there exist locally finite covers
$(U_n)_{n\in\N}$
and $(V_n)_{n\in \N}$
of~$M$ (resp., $N$)
by relatively compact, open sets
$U_n\sub M$ (resp., $V_n\sub N$)
such that $f(\gamma)|_{V_n}$ only depends
on $\gamma|_{U_n}$, i.e.,
\[
(\forall n\in \N)\;(\forall \gamma,\eta \in \Omega)\, \quad
\gamma|_{U_n}\,=\,\eta|_{U_n}\;\,\Rightarrow \,\;
f(\gamma)|_{V_n}\,=\, f(\eta)|_{V_n}\,.
\]
E.g., $f$ is almost local
if $M=N$ and $f$ is \emph{local}\index{local map}
in the sense that $f(\gamma)(x)$ only depends
on the germ of~$\gamma$ at $x\in M$.
As shown in \cite{DIF} (see also \cite[Theorem~10.4]{ZOO}),
almost locality
prevents pathologies as
in Example~\ref{extenso}.
\begin{proposition}\label{almloc}
Let $r,s,t\in \N_0\cup\{\infty\}$
and $f\colon C^r_c(M,E)\supseteq \Omega\to C^s_c(N,F)$
be an almost local map.
Assume that the restriction of~$f$ to
$\Omega\cap C^r_K(M,E)$
is $C^t$, for each compact set $K\sub M$.
Then $f$ is $C^t$.
\end{proposition}
An analogous result is available
for mappings between open subsets
of spaces of compactly supported sections
in vector bundles.
Almost local maps between
subsets of the space of compactly supported
smooth vector fields
occur in the construction of
the Lie group structure on $\Diff_c(M)$ (see \cite{DIF};
cf.\ \cite{DFR} and \cite{ZOO}).\\[2.5mm]
The \emph{proof of Proposition}~\ref{almloc}
exploits that the map
\[
\sigma \colon C^s(N,F)\to \bigoplus_{n\in \N} C^s(V_n,F)\,,\quad
\gamma\mto(\gamma|_{V_n})_{n\in \N}
\]
is a linear topological embedding with closed
image \cite[Proposition~8.13]{ZOO},
for each locally finite cover $(V_n)_{n\in \N}$
of~$N$ by relatively compact, open sets~$V_n$.
It hence suffices to show
that $\sigma \circ f$ is~$C^t$.
Let us assume that $\Omega=C^r_c(M,E)$ for simplicity.
There is a locally finite cover
$(\wt{U}_n)_{n\in \N}$ of~$M$ by relatively compact, open sets
such that $\wt{U}_n$ contains the closure of~$U_n$.
Let $h_n\colon \wt{U}_n \to \R$
be a compactly supported
smooth map such that $h_n|_{U_n}=1$.
Then the following map is $C^t$:
\[
f_n\colon C^r(\wt{U}_n,E)\to C^s(V_n,F)\, ,\quad
f_n(\gamma)\; :=\; f(h_n\cdot \gamma)|_{V_n}\,.
\]
Set $\rho \colon C^r_c(M,E)\to \bigoplus_{n\in\N}
C^r(\wt{U}_n,E)$, $\gamma\mto (\gamma|_{\wt{U}_n})_n$.
Then $\sigma\circ f=(\oplus_{n\in \N}f_n)\circ \rho$,
where $\oplus_{n\in \N}f_n$ is $C^t$ by Proposition~\ref{dsummp}.
Hence $\sigma\circ f$ and thus $f$ is~$C^t$.\smartqed\qed
\section{\,Lie group structures on directed unions
of Lie groups}\label{secconstr}
In some situations, it is possible
to construct Lie group structures
on ascending unions of Lie groups.\\[2.5mm]
{\bf Unions of finite-dimensional Lie groups.}
In the case of finite-dimensional manifolds,
an Extension Lemma
for Charts is available \cite[Lemma~2.1]{FUN}:

\emph{If $M$ and $N$ are finite-dimensional
$C^\infty$-manifolds such that $M\sub N$
and the inclusion map $M\to N$ is an
immersion, then each chart $\phi\colon U\to V$
of~$M$ which is defined on a relatively compact,
(smoothly) contractible subset $U\sub M$
extends to a chart of~$N$ on a domain
with analogous properties.}

Now consider a sequence $G_1\sub G_2\sub \cdots$
of finite-dimensional Lie groups
such that the inclusion maps
are smooth homomorphisms. Let $x\in G:=\bigcup_{n\in \N}G_n$,
say $x\in G_{n_0}$. We then pick a chart $\phi_{n_0}$
of $G_{n_0}$ around~$x$ whose domain
is relatively compact and contractible,
and use the extension lemma to obtain charts
$\phi_n$ of~$G_n$ for $n>n_0$
which are defined on relatively compact,
contractible
open sets, and such that $\phi_n$ extends $\phi_{n-1}$.
One then easily verifies (using Proposition~\ref{silvkom})
that the homeomorphisms $\phi:=\dl\, \phi_n$\vspace{-.7mm}
so obtained define a $C^\infty$-atlas
on $G$ (equipped with the direct limit topology),
which makes the latter a Lie group modelled
on $\dl\,\Lie(G_n)$\vspace{-.7mm}
(see \cite{FUN}).\footnote{Cf.\
\cite{NRW1}, \cite[Theorem~47.9]{KaM}
and \cite{DIR} for earlier,
less general results.}
By construction, $G=\bigcup_{n\in\N}G_n$
admits direct limit charts.
Moreover, it is clear from the construction
that
$G=\dl\,G_n$\vspace{-.7mm}
as a topological space and as a topological
group. Using Proposition~\ref{silvkom},
one easily infers that
$G=\dl\,G_n$\vspace{-.7mm}
also as a smooth manifold
and as a Lie group (see \cite[Theorem~4.3]{FUN}).
\begin{remark}
The preceding construction applies
just as well to ascending unions
of finite-dimensional smooth manifolds~$M_n$,
such that all inclusion maps
are immersions.\footnote{Compare \cite{Han}
for $\bigcup_{n\in \N}M_n$ as a topological
manifold.}
This enables
$G/H$ to be turned into the direct limit $C^\infty$-manifold
$\dl\,G_n/(H\cap G_n)$,\vspace{-.7mm}
for each closed subgroup $H\sub G$
(see \cite[Proposition~7.5]{FUN}).
Then the quotient map $G\to G/H$ makes $G$
a principal $H$-bundle over $G/H$,
using a suitable extension lemma
for sections in nested principal bundles
\cite[Lemma~6.1]{FUN}.

We mention that an equivariant version
of the above extension lemma (namely \cite[Lemma~1.13]{Wk07})
can be used to turn
the gauge group $\Gau(P)$
into a Lie group,
for each smooth principal bundle
$P\to K$
over a compact smooth manifold~$K$
whose structure group is a direct limit $G=\dl\,G_n$\vspace{-.7mm}
of finite-dimensional Lie groups (see
\cite[Lemma~1.14\,(e) and Theorem~1.11]{Wk07}).
\end{remark}
\begin{remark}\label{pathodl}
Direct limits $G=\dl\,G_n$\vspace{-.7mm}
of finite-dimensional Lie groups
are\linebreak
regular Lie groups in Milnor's sense \cite[Theorem~8.1]{FUN},
but they can be quite pathological in other ways.
E.g., the exponential map
$\exp_G=\dl\,\exp_{G_n}$\vspace{-.7mm}
need not be injective on any $0$-neighbourhood,
and
the exponential image
need not be an identity neighbourhood
in~$G$. Both pathologies occur
for
\[
G\; :=\; \C^{(\N)} \semid\; \R 
\; =\; \dl\,\C^n\semid\;\R\,,\vspace{-1mm}
\]
where $t\in \R$ acts on $\C^{(\N)}$
via $t.(z_k)_{k\in \N}:=(e^{ikt}z_k)_{k\in \N}$.
This can be checked
quite easily, using that
the exponential map of~$G$ is given explicitly by
$\exp_G((z_k)_{k\in \N},t)=\big(
\big(\frac{e^{ikt}-1}{ikt} z_k\big)_{k\in \N},t\big)$
(see \cite[Example~5.5]{DIR}).
\end{remark}
The preceding general construction
implies that every countably-dimensional
locally finite Lie algebra $\cg$
(i.e., each union $\cg=\bigcup_{n\in \N}\cg_n$
of finite-dimensional Lie algebras $\cg_1\sub
\cg_2\sub\cdots$),\index{locally finite Lie algebra}
when endowed with the finest locally convex vector
topology, arises as the Lie algebra
of some regular Lie group.\footnote{One chooses
a simply
connected Lie group $G_n$ with Lie algebra~$\cg_n$
and forms the direct limit group $G=\dl\,G_n$\vspace{-.7mm}
(see \cite[Theorem~5.1]{FUN} for the details).}
Such locally finite Lie algebras
have been much studied in recent years, e.g.
by\linebreak
Yu.\ Bahturin,
A.\,A. Baranov, G. Benkart, I. Dimitrov, K.-H. Neeb, I. Penkov, H. Strade,
N. Stumme, A.\,E. Zalesski\u{\i}, and others
(see \cite{BBZ04}, \cite{BB4}, \cite{DP04},
\cite{Ne00}, \cite{NS01}, \cite{PS03}, \cite{St99}
and the references therein).\\[3mm]
{\bf Unions of Banach--Lie groups.}
These are Lie groups under
additional\linebreak
hypotheses (which, e.g., exclude the pathologies described
in Remark~\ref{pathodl}).
\begin{theorem}\label{dahm2}
Let $G_1\sub G_2\sub\cdots$
be Banach--Lie groups
over $\K\in \{\R,\C\}$,
such that
all inclusion maps $\lambda_n \colon G_n\to G_{n+1}$
are $C^\infty_\K$-homomorphisms.
Set $G:=\bigcup_{n\in \N}G_n$.
Assume that {\rm(a)}--{\rm(c)}
are satisfied:
\begin{description}[(D)]
\item[\rm(a)]
For each $n\in \N$, there exists a norm $\|.\|_n$
on $\Lie(G_n)$ defining its topology,
such that $\|[x,y]\|_n \leq\|x\|_n\|y\|_n$
for all $x,y\in \Lie(G_n)$
and the continuous linear map
$\Lie(\lambda_n)\colon \Lie(G_n)\to \Lie(G_{n+1})$
has operator norm at most~$1$.
\item[\rm(b)]
The locally convex direct limit topology on
$\cg:=\bigcup_{n\in \N}\Lie(G_n)$ is Hausdorff.
\item[\rm(c)]
$\exp_G:=\dl\, \exp_{G_n}\colon \cg\to G$\vspace{-.7mm}
is injective on some $0$-neighbourhood.
\end{description}
Then there exists a $\K$-analytic Lie group
structure on~$G$ which makes $\exp_G$
a $\K$-analytic local diffeomorphism at~$0$.
If, furthermore, $\cg=\bigcup_{n\in \N}\Lie(G_n)$ is
compactly regular, then $G$
is a regular Lie group in Milnor's sense.
\end{theorem}
\emph{Sketch of proof.}
The Lie group structure is constructed in \cite{Dah},
along the following lines:
Applying Dahmen's Theorem~\ref{thmdah}
(to the complexification $\cg_\C$, if $\K=\R$),
one finds that the Baker--Campbell--Hausdorff (BCH-)
series converges to a $\K$-analytic
map
$\bigcup_{n\in \N}B^{\Lie(G_n)}_r(0)\times B^{\Lie(G_n)}_r(0)\to \cg$
for some $r>0$.
Because $\exp_G$ is locally injective,
it induces an isomorphism~$\phi$ of local
groups from some $0$-neighbourhood
$U\sub \bigcup_{n\in \N}B^{\Lie(G_n)}_r(0)$
onto some subset~$V$ of~$G$.
We give $V$ the $C^\omega_\K$-manifold
structure making $\phi$ a $C^\omega_\K$-diffeomorphism.
Now standard arguments can be used to make~$G$
a Lie group with~$V$ as an open
submanifold.
The proof of regularity will be sketched
in Section~\ref{secregu}.\vspace{2mm}\smartqed\qed

\noindent
The author does not know whether the Lie groups $G$
in Theorem~\ref{dahm2}
are\linebreak
always the direct limit $\dl\, G_n$\vspace{-.7mm}
in the category of Lie groups
(unless additional hypotheses are satisfied).\\[3mm]
{\bf Another construction principle.}
There is another
construction principle for a Lie group
structure on a union $G=\bigcup_{n\in \N}G_n$
of Lie groups (or a group which is a union
$G=\bigcup_{n\in \N}M_n$ of manifolds),
which produces
Lie groups modelled on Silva spaces
or ascending unions of $k_\omega$-spaces.
A direct limit Lie group structure
can be constructed on~$G$ if (1) there
are compatible charts $\phi_n$
of the Lie groups $G_n$ (resp., the manifolds~$M_n$)
around each point in~$G$;
and (2) suitable
hypotheses are satisfied
which ensure that the transition maps
between charts of the form $\dl\,\phi_n$\vspace{-.7mm}
are $C^\infty_\K$,
because they are mappings of the form discussed
in Proposition~\ref{silvkom} (see \cite[Lemma~14.5]{COM}).
\section{\,Examples of directed unions of
Lie groups}\label{secuninf}
The main examples
of ascending unions of infinite-dimensional
Lie groups were already briefly described
in the introduction.
We now provide more\linebreak
details.
Notably, we discuss
the existence of direct limit charts,
and compact regularity.
As already mentioned,
the latter gives information
on the homotopy groups (see (\ref{dlhomot}))
and can help to verify regularity
in Milnor's sense (see Theorem~\ref{dahm2}
and Section~\ref{secregu}).
A special case of
\cite[Corollary~3.6]{HOM} is useful.
\begin{lemma}\label{thsreg}
If the Lie group
$G=\bigcup_{n\in \N}G_n$
admits a weak direct limit chart,
then $G=\bigcup_{n\in \N}G_n$
is compactly regular if and only if
$\Lie(G)=\bigcup_{n\in \N}\Lie(G_n)$
is compactly regular.
\end{lemma}
In the case of
an (LF)-space $E=\bigcup_{n\in\N} E_n$,
there is a quite concrete characterization of compact regularity
only in terms of properties of the steps~$E_n$
(see \cite[Theorem~6.4 and its corollary]{Wen}):\index{Wengenroth's Theorem}
\begin{theorem}\label{wengen}
Let $E_1\sub E_2\sub \cdots$ be Fr\'{e}chet
spaces, with continuous linear inclusion maps.
Give $E=\bigcup_{n\in \N}E_n$
the locally convex direct limit topology.
Then $E=\bigcup_{n\in \N}E_n$
is compactly regular
if and only if
for each $n\in \N$, there exists $m\geq n$
such that for all $k\geq m$,
there is a $0$-neighbourhood $U$ in $E_n$
on which $E_k$ and $E_m$ induce the
same topology.
In this case, $E$ is also boundedly regular and complete.
\end{theorem}
We mention that a Hausdorff (LF)-space
is boundedly regular if and only if it is
Mackey complete \cite[1.4\,(f), p.\,209]{Fl80}.\\[3mm]
{\bf Groups of compactly supported
diffeomorphisms.}
The Lie group\linebreak
$\Diff_c(M)=\bigcup_{n\in \N}\Diff_{K_n}(M)$
(discussed in the introduction)
admits a direct limit chart (cf.\ \cite[\S5.1]{COM}).
Moreover, the LF-space $\cV_c(M)=\bigcup_{n\in \N}\cV_{K_n}(M)$
is compactly regular (see \ref{baslcx}\,(c))
and hence also $\Diff_c(M)$
(by Lemma~\ref{thsreg}).
To avoid exceptional cases
in our later discussions
of direct limit properties,
we assume henceforth that $M$ is non-compact
and of positive dimension.\\[2.5mm]
{\bf Test function groups.}
Let $M$ and an exhaustion $K_1\sub K_2\sub\cdots$ of~$M$
be as in the definition of $\Diff_c(M)$,
$H$ be a Lie group modelled
on a locally convex space,
and $r\in \N_0\cup\{\infty\}$.
We consider the ``test function group''
$C^r_c(M,H)$ of
$C^r$-maps $\gamma\colon M\to H$
such that the closure of $\{x\in M\colon \gamma(x)\not={\bf 1}\}$
(the support of $\gamma$) is compact.
Let $C^r_{K_n}(M,H)$
be the subgroup of functions supported in~$K_n$.
Then $C^r_{K_n}(M,H)$ is a Lie group modelled
on $C^r_{K_n}(M,\Lie(H))$, and
$C^r_c(M,H)$ is a Lie group
modelled on the locally convex direct
limit $C^r_c(M,\Lie(H))=\dl\,C^r_{K_n}(M,\Lie(H))$\vspace{-.7mm}
(\cite{GCX}; cf.\ \ref{introhis}
for special cases). Also,
\[
C^r_c(M,H)\;=\; {\textstyle \bigcup_{n\in \N}}\, C^r_{K_n}(M,H)
\]
admits a direct limit chart (cf.\ \cite[\S7.1]{COM}).
Furthermore, $C^r_c(M,\Lie(H))\!=$ $\bigcup_n\, C^r_{K_n}(M,\Lie(H))$
is compactly regular as a consequence
of \ref{baslcx}\,(c).
We now assume
that $H$ is non-discrete
and $M$ non-compact,
of positive dimension.\\[2.5mm]
{\bf Weak direct products of Lie groups.}
Given a sequence
$(H_n)_{n\in \N}$ of Lie groups,
its weak direct product
$G:=\prod_{n\in \N}^*H_n$
(as introduced before Proposition~\ref{propcriter})
has a natural Lie group structure~\cite[\S7]{MEA},
modelled on the locally convex direct
sum $\bigoplus_{n\in\N}\Lie(H_n)$.
Then $G=\bigcup_{n\in \N}G_n$, identifying the partial
product $G_n:=\prod_{k=1}^nH_k$ with a subgroup
of~$G$. By construction, $G=\bigcup_{n\in \N}G_n$
has a direct limit chart. Furthermore,
$\Lie(G)=\bigoplus_{n\in \N}\Lie(H_n)
=\dl\, \Lie(G_n)$\vspace{-.5mm}
is compactly regular,
as locally convex direct sums are
boundedly regular \cite[Ch.\,3, \S1, no.\,4, Proposition~5]{BTV}
and induce the given topology on each finite
partial product
(cf.\ Propositions~7 or 8\,(i)
in \cite[Ch.\,2, \S4, no.\,5]{BTV}).
To avoid exceptional cases,
we assume henceforth
that each $H_n$ is non-discrete.\\[2.5mm]
{\bf Unit groups of
unions of Banach algebras.}
Let
$A_1\sub A_2\sub\cdots$ be
unital complex Banach algebras
(such that all inclusion maps
are continuous homomorphisms of unital algebras).
Give
$A:=\bigcup_{n\in \N}A_n$
the locally convex direct limit topology.
Then $A^\times$ is open in~$A$
and if $A$ is Hausdorff (which we assume now),
then
$A^\times$ is a complex Lie group
\cite[Proposition~12.1]{COM}.
Moreover, $A^\times=\bigcup_{n\in \N}A_n^\times$,
and the identity map $\id_{A^\times}$ is a direct
limit chart (cf.\ \cite{DaW} and \cite{Eda}
for related results).

If each inclusion map $A_n\to A_{n+1}$
is a topological embedding or each inclusion map
a compact
operator, then $A=\bigcup_{n\in \N}A_n$
and hence also $A^\times=\bigcup_{n\in \N}A_n^\times$
is compactly regular.
However, for particular
choices of the steps, $A=\bigcup_{n\in \N}A_n$
is not compactly regular
(see \cite[Example~7.8]{HOM},
based on \cite[Remark~1.5]{BMS}).\\[2.5mm]
{\bf Lie groups of germs of analytic mappings.}
Let $H$ be a
complex Banach-Lie group,
$\|.\|$ be a norm on $\Lie(H)$ defining its topology,
$X$ be a complex
metrizable locally convex space
and $K\sub X$ be a non-empty compact set.
Let $W_1\supseteq W_2\supseteq \cdots$
be a fundamental sequence of open neighbourhoods of~$K$
in~$X$ such that each connected
component of $W_n$ meets~$K$.
Then the set
$\Germ(K,H)$
of germs around~$K$ of $H$-valued complex analytic functions
on open neighbourhoods of~$K$ can be made a Lie group
modelled on the locally convex direct limit
\[
\Germ(K,\Lie(H))\; =\; \dl\, \Hol_b(W_n,\Lie(H))
\]
of the Banach spaces
$\cg_n:=\Hol_b(W_n,\Lie(H))$ of bounded $\Lie(H)$-valued
complex analytic functions on~$W_n$,
equipped with the supremum norm (see~\cite{HOL}).
The group operation arises from
pointwise multiplication of representatives of germs.
The identity component $\Germ(K,H)_0$
is the union
\[
\Germ(K,H)_0\;=\; \bigcup_{n\in \N}G_n
\]
of the Banach--Lie
groups $G_n:=\langle [\exp_H\circ\, \gamma]\colon \gamma\in \cg_n\rangle$,
and $\Germ(K,H)_0=\bigcup_{n\in\N}G_n$
admits a direct limit chart~\cite[\S10.4]{COM}.
Theorem~\ref{wengen}
implies that $\Germ(K,\Lie(H))=\bigcup_{n\in \N}\cg_n$
is compactly regular
(see \cite{REG}),
and thus $\Germ(K,H)_0=\bigcup_{n\in \N}G_n$
is compactly regular (see already
\cite[Theorems~21.15 and 21.23]{Ch85}
for the bounded regularity and completeness
of $\Germ(K,\Lie(H))$ if $X$ is a normed
space; cf.\ \cite{Mj79}).
In the most relevant case
where $X$ and~$H$ are finite-dimensional,
we can choose $W_{n+1}$ relatively compact in~$W_n$.
Then the restriction maps
$\Hol_b(W_n,\Lie(H))\to \Hol_b(W_{n+1},\Lie(H))$
are compact
operators~\cite[\S10.5]{COM}
and thus $\Germ(K,\Lie(H))$
is a Silva space.\\[2.5mm]
{\bf Lie groups of germs of analytic diffeomorphisms.}
If $X$ is a complex
Banach space
and $K\sub X$ a non-empty compact set,
let $\GermDiff(K,X)$
be the set of germs around~$K$ of
$\C$-analytic diffeo\-morphisms
$\gamma\colon U\to V$
between open neighbourhoods $U$ and $V$ of~$K$
(which may depend on $\gamma$),
such that $\gamma|_K=\id_K$.
Then $\GermDiff(K,X)$
is a Lie group
modelled on the locally convex direct limit
\[
\Germ(K,X)_K\; :=\; \dl\, \Hol_b(W_n,X)_K\,,
\]
where $W_n$ and $\Hol_b(W_n,X)$
are as in the last example
and $\Hol_b(W_n,X)_K:=\{\zeta \in \Hol_b(W_n,X)\colon \zeta|_K=0\}$
(see \cite[\S15]{COM} for the case
$\dim(X)<\infty$, and \cite{Dah} for the general
result).
The group operation arises from
composition of representatives of germs.
Now the set $M_n$
of all elements of $\GermDiff(K,X)$
having a representative in $\Hol_b(W_n,X)_K$
is a Banach manifold, and
\[
\GermDiff(K,X)\;=\; \bigcup_{n\in \N}M_n
\]
has a direct limit chart (see \cite{Dah}; cf.\
\cite[Lemma~14.5 and \S15]{COM}).\linebreak
$\GermDiff(K,X)\!=\!\bigcup_n \!M_n$
is compactly regular by Theorem~\ref{wengen}
and Lemma~\ref{thsreg} (see~\cite{Dah});
if $X$ is finite-dimensional,
then $\Germ(K,X)_K$ is a Silva space.\\[2.5mm]
{\bf Unions of Lie groups modelled on Sobolev spaces.}
The Lie groups
$H^{\downarrow s}(K,F)
=\bigcup_{n\in \N}H^{s+\frac{1}{n}}(K,F)$
(as in the introduction)
are studied in the work in progress~\cite{REG}.
By construction, they
admit a direct limit chart,
and they are modelled
on the Silva space
$H^{\downarrow s}(K,\Lie(F))
=\bigcup_{n\in \N}H^{s+\frac{1}{n}}(K,\Lie(F))$
(and hence compactly regular).
We mention that the Lie group structure
on $H^{\downarrow s}(K,F)$
can be obtained via Theorem~\ref{dahm2};
therefore
$H^{\downarrow s}(K,F)$ is a\linebreak
regular
Lie group in Milnor's sense.
Compare \cite{Pi08} (in this volume) for\linebreak
analysis
and probability theory
on variants of the Lie groups $H^s(K,F)$
(with $s>\dim(K)/2$), and limit processes
as $s\downarrow \dim(K)/2$.
\section{\,Direct limit properties of ascending unions}\label{secdlprop}
We now discuss the direct limit properties
of ascending unions of infinite-dimensional
Lie groups in the categories
of Lie groups, topological groups,
smooth manifolds and topological spaces.\vfill\pagebreak

\noindent
{\bf Tools to prove or disprove direct
limit properties.}
Such tools were
provided in~\cite{COM}.
Recall that a real locally convex space~$E$
is said to be \emph{smoothly regular}\index{smoothly regular space}
(or: to \emph{admit smooth bump
functions})\index{smooth bump function}\index{bump function}
if the topology on~$E$
is initial with respect to $C^\infty(E,\R)$.
\begin{remark}
If $U\sub E$
is a $0$-neighbourhood and the topology is initial with respect to
$C^\infty(E,\R)$, then $\bigcap_{j=1}^n
f_j^{-1}(]{-\ve},\ve[) \sub U$
for suitable $\ve>0$
and $f_1,\ldots, f_n\in C^\infty(E,\R)$
such that $f_1(0)=\cdots=f_n(0)=0$.
Then $f^{-1}(]{-\delta},\delta[)\sub U$
with $f:=f_1^2+\cdots+ f_n^2$ and $\delta:=\ve^2$.
Let $g\colon \R\to \R$
be a smooth function such that $g(\R)\sub [0,1]$,
$g(0)=1$
and $g(x)=0$ if $|x|\geq\delta/2$.
Then $h:=g\circ f\colon E\to \R$
is a smooth function such that $h(0)=1$
and $\Supp(h)\sub U$
(a ``smooth bump function''
supported in~$U$). This explains
the terminology.
\end{remark}
\begin{example}
Every Hilbert space~$H$
admits smooth bump functions
(because $H\to \R$, $x\mto \|x\|^2$
is smooth).
As a consequence,
every locally convex space
which admits a linear topological embedding
into a direct product of Hilbert spaces
(for example, every nuclear locally convex space)
admits smooth bump functions
(cf.\ also \cite[Chapter~III]{KaM}).
\end{example}
\begin{proposition}\label{tools2}
Consider a Lie group $G=\bigcup_{n\in \N}G_n$,
where $G_1\sub G_2\sub\cdots$
are Lie groups and all inclusion maps
$G_n\to G_{n+1}$ and $G_n\to G$
are smooth homomorphisms.
Assume that $G=\bigcup_{n\in \N}G_n$
admits a direct limit chart.
Then the following hold:
\begin{description}[(D)]
\item[\rm(a)]
If $G=\dl\,G_n$
as a topological group,
then $G=\dl\,G_n$\vspace{-.7mm}
as a Lie group.
\item[\rm(b)]
$G=\dl\,G_n$
as a topological space
if and only if $\Lie(G)=\dl\,\Lie(G_n)$\vspace{-.7mm}
as a topological space.
\item[\rm(c)]
If $\Lie(G)$ admits smooth bump functions,
then $G=\dl\,G_n$\vspace{-.7mm}
as a $C^\infty_\R$-manifold
if and only if $\Lie(G)=\dl\,\Lie(G_n)$\vspace{-.7mm}
as a $C^\infty_\R$-manifold.\vspace{1mm}
\end{description}
\end{proposition}
{\bf Direct limit properties
of the main examples.}
Using Proposition~\ref{tools2},
Proposition~\ref{propcriter}
(to recognize direct
limits of topological groups)
and a\linebreak
counterpart of
Proposition~\ref{silvkom}
for analogous ascending unions
of manifolds~\cite[Proposition~9.8]{COM},
one obtains the following information
concerning the direct limit properties
of the examples from
in Section~\ref{secuninf} (see \cite{COM};
the
properties of $H^{\downarrow s}(K,F)$
follow from \cite[Proposition~9.8]{COM}).

The
entries
in the following table
indicate whether $G=\dl\,G_n$\vspace{-.7mm}
holds in the category shown on the left,
for the Lie group described at the top.
The abbreviation ``dep''
is used if the answer depends on
special properties
of the group(s) involved.
We abbreviate
``category'' by ``cat,''
``group'' by ``gp,''
``space'' by ``sp,''
``topological'' by ``top'',\index{direct limit property}
and ``smooth manifold'' by ``mfd.''\vfill\pagebreak

\noindent\hspace*{-1.7mm}
{\footnotesize
\begin{tabular}{||c||c|c|c|c|c|c|c||}\hline\hline
cat$\backslash$gp & $\text{Diff}_c(M)$ 
& $C^\infty_c(\hspace*{-.2mm}M\hspace*{-.2mm},\hspace*{-.5mm}H)$
& $\prod^*_n \hspace*{-.7mm}H_n$ &
$A^\times$ & $\Germ(K\hspace*{-.2mm},\hspace*{-.4mm}H)_0 $ &
$\GermDiff(\hspace*{-.2mm}K,\hspace*{-.4mm}X)$
& $H^{\downarrow s}(K,\hspace*{-.4mm}F)$ \\ \hline\hline
Lie gps& yes  & yes & yes & yes & yes & ---\, & yes \\ \hline
top gps & yes   & yes & yes & yes & yes & ---\, & yes \\ \hline
mfds & no & no & dep* & dep**& \,yes$^\dag$
& \,yes$^{\dag\dag}$& yes \\ \hline
top sps & no & no & dep* & dep** & \,yes$^{\dag}$ & \,yes$^{\dag\dag}$
& yes \\ \hline\hline
\end{tabular}}\\[4.6mm]
{\footnotesize
\hspace*{1mm}* \hspace*{.5mm}``yes'' if each $H_n$ is finite-dimensional
or modelled on a $k_\omega$-space;
``no'' if each~$H_n$\linebreak
\hspace*{5mm}is modelled
on an infinite-dimensional Fr\'{e}chet space
(which we assume nuclear\linebreak
\hspace*{5mm}when dealing with the category
of smooth manifolds). Other cases unclear.\\[1mm]
** \hspace*{-.8mm}``yes''\hspace*{-.3mm} if each $A_n$ is finite-dimensional
or each inclusion map $\lambda_n\colon A_n\to A_{n+1}$ a\linebreak
\hspace*{5mm}compact operator;
``no'' (when dealing with the category of topological
spaces), if\linebreak
\hspace*{5mm}$A_n$ is infinite-dimensional,
$A_n\subset A_{n+1}$ and $\lambda_n$ a topological
embedding for\linebreak
\hspace*{5mm}each~$n$. Other cases unclear.\\[1mm]
\hspace*{.8mm}$\dag$ \hspace*{.5mm}``yes'' if $X$ and $H$
are finite-dimensional;
general case unknown.\\[1mm]
${\dag\dag}$ \hspace*{.2mm}``yes'' if $X$ is finite-dimensional;
general case unknown.}
\section{\,Regularity in Milnor's sense}\label{secregu}
Experience tells that
if one tries to prove regularity in Milnor's
sense for a Lie group
$G=\bigcup_{n\in \N}G_n$,
then regularity of the Lie groups~$G_n$
does not suffice to carry out the desired
arguments. But strengthened
regularity properties
increase the chances for success.
\begin{definition}
{\rm Given $k\in \N_0$, we say that
a Lie group~$G$ is \emph{$C^k$-regular}
if it is a regular Lie group
in Milnor's sense and
\[
\evol_G\colon C^\infty([0,1],\Lie(G))\to G
\]
is smooth with respect to the
$C^k$-topology on $C^\infty([0,1],\Lie(G))$
(induced by $C^k([0,1],\Lie(G))$).
If each $\gamma\in C^k([0,1],\Lie(G))$
has a product integral $\eta_\gamma$
and the map
\[
\evol_G\colon C^k([0,1],\Lie(G))\to G\,,\quad
\gamma\mto \eta_\gamma(1)
\]
is smooth,
then we say that the Lie group~$G$
is \emph{strongly
$C^k$-regular}.}\index{$C^k$-regular Lie group}\index{strongly
$C^k$-regular Lie group}
\end{definition}
E.g., every Banach--Lie group
is strongly $C^0$-regular~\cite{GaN}.
Although much of the following
remains valid for $C^k$-regular Lie groups,
we shall presume
strong $C^k$-regularity,
as this simplifies the presentation.
We also suppress possible variants involving
bounded regularity instead
of compact regularity.
All results presented in this section
are taken from~\cite{REG}.

In the regularity proofs for our main classes
of direct limit groups, we always use an isomorphism
$C^k([0,1],\dl\,E_n)\isom \dl\,C^k([0,1],E_n)$\vspace{-.7mm}
at a pivotal point.
Let us begin with the elementary
case of locally convex direct sums.
\begin{lemma}\label{dsanddl}
If $(E_n)_{n\in \N}$ is a sequence
of Hausdorff locally convex spaces,
then
$C^k([0,1],\bigoplus_{n\in \N}E_n)=\bigoplus_{n\in \N}C^k([0,1],E_n)$,
for all $k\in \N_0$.
\end{lemma}
\emph{Sketch of proof.}
The locally convex direct sum $\bigoplus_{n\in \N}E_n
=\bigcup_{n\in \N} (E_1\times\cdots\times E_n)$
is compactly regular, because it is boundedly regular
by \cite[Chapter~3, \S1, no.\,4, Proposition~5]{BTV}
and induces the given topology on each finite partial
product (cf.\ Propositions~7 or~8\,(i)
in \cite[Chapter~2, \S4, no.\,5]{BTV}).
Therefore\linebreak
$C^k([0,1],\bigoplus_{n\in \N}E_n)$
and $\bigoplus_{n\in \N}C^k([0,1],E_n)$
coincide as sets.
Comparing\linebreak
$0$-neighbourhoods,
we see that both vector topologies coincide
(using that boxes are typical $0$-neighbourhoods in a countable
direct sum).\smartqed\qed
\begin{remark}
Although
$C^\infty([0,1],\R^{(\N)})=\bigcup_{n\in \N}C^\infty([0,1],\R^n)$
as a set,
the topology on the left hand side
is properly coarser than the locally convex
direct limit topology $\cO_{\lcx}$ on the right hand side,
because
\[
\Big\{\gamma=(\gamma_n)_{n\in \N}\in C^\infty([0,1],\R^{(\N)})\colon
(\forall n\in \N) \; {\textstyle \frac{d^n\gamma_n}{dx^n}}(0)\in \;]{-1,1}[\,\Big\}\;
\in\; \cO_{\lcx}
\]
is not a $0$-neighbourhood
in $C^\infty([0,1],\R^{(\N)})
=\pl_{m\in \N_0}\, C^m([0,1],\R^{(\N)})$.\vspace{-.7mm}
Thus Lemma~\ref{dsanddl}
becomes false for $k=\infty$,
explaining the need for
$C^k$-regularity with
finite~$k$.
\end{remark}
{\bf Weak direct products of Lie groups.}
If $k\in \N_0$ and $(H_n)_{n\in \N}$ is a sequence
of strongly $C^k$-regular Lie groups,
then $\prod_{n\in\N}^*H_n$ is strongly $C^k$-regular
(and hence regular) since its evolution map
can be obtained as the composition\vspace{-1mm}
\[
C^k\big([0,1],{\textstyle \bigoplus_n \Lie(H_n)}\big)\,
\stackrel{\isom}{\longrightarrow}\,
{\textstyle \bigoplus_n C^k([0,1], \Lie(H_n))}
\, \stackrel{\oplus_{n\in \N}\evol_{H_n}}{\longrightarrow}\,
{\textstyle \prod_n^*}\, H_n
\]
(cf.\ Proposition~\ref{dsummp} for the definition
and smoothness of $\oplus_n \evol_{H_n}$).\\[2.5mm]
{\bf Test function groups.}
Given a $\sigma$-compact, finite-dimensional
smooth\linebreak
manifold~$M$ and a $C^k$-regular
Lie group~$H$, pick
a locally finite family $(M_n)_{n\in \N}$ of compact
submanifolds with boundary of~$M$,
the interiors of which cover~$M$.
Then standard arguments (based
on suitable exponential laws for function spaces)
show that $H_n:=C^r(M_n,H)$ is $C^k$-regular,
for each $n\in \N$.
The map $\sigma\colon C^r_c(M,\Lie(H))\to \bigoplus_n C^r(M_n,\Lie(H))$,
$\gamma\mto (\gamma|_{M_n})_{n\in\N}$
is continuous linear and hence also
the map $\tau:= C^k([0,1],\sigma)$ from
$C^k([0,1],C^r_c(M,\Lie(H)))$ to
$C^k([0,1],\bigoplus_n C^r(M_n,\Lie(H)))
\isom \bigoplus_n C^k([0,1], C^r(M_n,\Lie(H)))$.
Furthermore,
\[
\rho \colon G:=C^r_c(M,H)\to {\textstyle \prod^*_{n\in \N}}C^r(M_n,H)\,,
\quad \gamma\mto (\gamma|_{M_n})_{n\in \N}
\]
is an isomorphism of Lie groups onto a
closed Lie subgroup (and embedded submanifold)
of the weak direct product~$P$.
Using point evaluations, one finds that
the composition
\[
\evol_P\circ \, \tau \, =\, \oplus_n \evol_{H_n}\circ \, \tau
\colon C^r_c(M,H)\to {\textstyle \prod^*_{n\in \N}}C^r(M_n,H)
\]
(which is smooth by the
preceding example)
takes its image in the image of $\rho$.
Then $f:=\rho^{-1}\circ \evol_P\circ \, \tau\colon C^k([0,1],C^r_c(M,\Lie(H)))
\to C^r_c(M,H)$ is a smooth map,
and one verifies using point evaluations
that $f=\evol_G$.\\[2.5mm]
A similar (but more complicated)
argument shows that $\Diff_c(M)$ is regular.\\[2.5mm]
{\bf Ascending unions of Banach--Lie groups.}
As a preliminary,
observe that regularity
in Milnor's sense (and strong $C^k$-regularity)
can be defined just as well for \emph{local}
Lie groups~$G$; in this case, one requires
that a smooth\linebreak
evolution $\evol_G$
exists on some open $0$-neighbourhood in $C^\infty([0,1],\Lie(G))$
(resp., $C^k([0,1],\Lie(G))$). In the case
of global Lie groups,
the local notions of regularity are equivalent
to the corresponding global ones (see \cite{REG}
and \cite{GaN}; cf.\ \cite[lemma on p.\,409]{KaM}).
See \cite{Sme} for
the next theorem
(and \cite{Muj} or \cite[Theorem~I.7.2]{Sm83}
for a variant beyond compact regularity).
\begin{theorem}\label{Muji}
Consider a Hausdorff locally convex space~$E$
which is the\linebreak
locally convex direct limit
of Hausdorff locally convex spaces
$E_1\sub E_2\sub \cdots$. If $E=\bigcup_{n\in \N}E_n$
is compactly regular,
then the natural continuous linear map
\[
\dl\, C([0,1],E_n)\to C([0,1],E)
\]
is an isomorphism of topological vector spaces.
\end{theorem}
\begin{remark}\label{bonet}
If $E$ is a locally convex space
which is \emph{integral complete},\footnote{That is,
every continuous curve $\gamma\colon [0,1]\to E$
has a Riemann integral in~$E\,$~\cite{LS00}.}
then
\begin{equation}\label{neednba}
C^k([0,1],E)\isom E^k\times C([0,1],E)
\end{equation}
naturally via
$\gamma\mto (\gamma(0),\ldots, \gamma^{(k-1)}(0),\gamma^{(k)})$,
for each $k\in \N$.
If $E=\dl\,E_n$\vspace{-.7mm}
is a locally convex direct limit
of integral complete locally convex spaces
and $E=\bigcup_{n\in\N}E_n$ is compactly regular,
then also~$E$ is integral complete.
Moreover,~(\ref{neednba}),
\ref{baslcx}\,(e)
and Theorem~\ref{Muji} imply
that
\begin{equation}\label{needaga}
\dl\,C^k([0,1],E_n)\, \isom \, C^k([0,1],\dl\,E_n)\,.\vspace{-.7mm}
\end{equation}
Alternatively,
(\ref{needaga}) follows from Theorem~\ref{Muji}
because $C^k([0,1],E)\isom C([0,1],E)$
naturally if $E$ is integral complete~\cite{REG},
as follows from (\ref{neednba})
and the fact that $C([0,1],E)\isom E^m\times C([0,1], E)$
naturally for each $m\in \N$,
by a suitable (elementary) variant of Miljutin's Theorem~\cite{Miu}
provided in~\cite{REG}.
\end{remark}
\begin{numba}\label{nowprove}
Assume that
$\cg=\bigcup_{n\in \N}\Lie(G_n)$ is compactly regular
in Theorem~\ref{dahm2}.\linebreak
Following~\cite{REG},
we now explain in the essential case where $\K=\C$
that $G$ is regular in Milnor's sense.\footnote{The real case
follows easily via complexification
on the level of local groups.}
Let $r>0$ be as in the earlier parts
of the proof and $U_n:=B^{\Lie(G_n)}_r(0)$.
Because the BCH-series
has the same shape for each $n\in\N$,
one finds $s>0$
such that an evolution
$\evol_{U_n}$ exists as a map from
$C([0,1], B^{\Lie(G_n)}_s(0))$ to $U_n$,
for each $n\in\N$. Since~$U_n$ is bounded,
Theorem~\ref{thmdah} shows that
$\evol_U:=\dl\,\evol_{U_n}\colon \bigcup_{n\in \N} C([0,1], B^{\Lie(G_n)}_s(0))
\to U:=\bigcup_{n\in \N}U_n$ is $\C$-analytic
and hence also $\exp_G\circ \evol_U$,
which is a local group version of the
evolution map for~$G$. Hence $G$ is regular and in fact
strongly $C^0$-regular.\smartqed\qed
\end{numba}
Using Theorem~\ref{dahm2},
one readily deduces that
$\Germ(K,H)$ and $H^{\downarrow s}(K,F)$
are strongly $C^0$-regular,
and also $A^\times=\bigcup_{n\in \N}A_n^\times$
if $A=\bigcup_{n\in \N}A_n$ is
compactly regular (see \cite{REG}).
The proof of compact regularity
for $\GermDiff(K,X)$ is more involved,
but eventually also boils down to
Theorem~\ref{thmdah}
(see \cite{Dah}).\\[3mm]
{\bf An idea which might lead to non-regular Lie groups.}
An observation from~\cite{REG} might be a source
of Lie groups which are not regular in Milnor's
sense
although they are modelled
on Mackey complete locally convex spaces:\index{non-regular Lie group}
\begin{proposition}\label{pathrg}
Suppose that, for each $n\in \N$,
there exists a Lie group~$H_n$
modelled on a Mackey complete
locally convex space
which is regular but
not $C^n$-regular
because
$\evol_{H_n}\colon C^\infty([0,1],\Lie(H_n))\to H_n$
is discontinuous with respect to the
$C^n$-topology.
Then $G:=\prod_{n\in \N}^*H_n$ is a Lie group
modelled on the Mackey complete
locally convex space $\Lie(G)=\bigoplus_{n\in \N}\Lie(H_n)$.
It has
an evolution $\evol_G\colon C^\infty([0,1],\Lie(G))\to G$,
but $\evol_G$ fails to be continuous
and thus $G$ is not a regular Lie group
in Milnor's sense.
\end{proposition}
\section{\,Homotopy groups of ascending unions
of Lie groups}\label{sechomotop}
We have seen that all
main examples of ascending unions $G=\bigcup_{n\in \N}G_n$
of Lie groups admit a direct limit chart,
and thus
\begin{equation}\label{again}
\pi_k(G)\;=\; \dl\, \pi_k(G_n)\vspace{-.7mm}\quad\mbox{for all $\,k \in \N_0$}
\end{equation}
(see \ref{subsechom}).
Alternatively, many (but not all)
of them are compactly regular.
In this case, (\ref{again})
holds by an elementary argument,
but one has to pay the price
that the proof of compact regularity may require
specialized functional-analytic tools
(like Wengenroth's theorem
recalled above).

It is an interesting feature that
the approach via (weak) direct limit charts
even extends to Lie groups~$G$
in which an ascending union $\bigcup_{n\in \N}G_n$
is merely \emph{dense} (and to similar,
more general situations). Weak direct
limit charts (as defined in \ref{susecdlcha})
have to be replaced by certain ``well-filled
charts'' then. The precise setting
will be described now.
Besides smooth manifolds,
it applies to topological manifolds
and more general
topological spaces (like
manifolds with boundary or corners).
Given a subset~$A$ of a real vector space~$V$,
let us write $\conv_2(A):=\{tx+(1-t)y\colon x,y\in A, t\in [0,1]\}$.
\begin{definition}
{\rm Let $M$ be a topological space
and $(M_\alpha)_{\alpha\in A}$ be a directed family of
topological spaces such that
$M_\infty:=\bigcup_{\alpha\in A}M_\alpha$
is dense in~$M$ and
all inclusion maps $M_\alpha\to M$ and
$M_\alpha\to M_\beta$ (for $\alpha\leq\beta$)
are continuous.
We say that a homeomorphism
$\phi\colon U\to V\sub E$
from an open subset
$U \sub M$ onto an arbitrary subset~$V$
of a topological vector space~$E$
is a \emph{well-filled chart} of~$M$\index{well-filled chart}
if there exist
$\alpha_0\in A$
and homeomorphisms
$\phi_\alpha\colon U_\alpha\to V_\alpha\sub E_\alpha$
from open subsets $U_\alpha\sub M_\alpha$
onto subsets $V_\alpha$ of certain topological
vector spaces $E_\alpha$ for $\alpha\geq\alpha_0$
such that the following conditions
are satisfied:
\begin{description}[(D)]
\item[(a)]
$E_\alpha\sub E$,
$E_\alpha\sub E_\beta$
if $\alpha\leq \beta$
and the inclusion maps $E_\alpha\to E$
and $E_\alpha\to E_\beta$ are continuous
and linear.
\item[(b)]
For all $\alpha\geq\alpha_0$, we have $U_\alpha\sub U$
and $\phi|_{U_\alpha}=\phi_\alpha$.
\item[(c)]
For all $\beta\geq \alpha\geq\alpha_0$, we have
$U_\alpha\sub U_\beta$
and $\phi_\beta|_{U_\alpha}=\phi_\alpha$.
\item[(d)]
$U_\infty:=\bigcup_{\alpha\geq\alpha_0}U_\alpha=U\cap M_\infty$.
\item[(e)]
There exists a non-empty (relatively) open set $V^{(2)}\sub V$
such that $\conv_2(V^{(2)})\sub V$
and $\conv_2(V_\infty^{(2)})\sub V_\infty$,
where $V_\infty:=\bigcup_{\alpha\geq\alpha_0}V_\alpha$
and $V_\infty^{(2)}:=V^{(2)}\cap V_\infty$.
\item[(f)]
For each $\alpha\geq\alpha_0$ and compact set
$K\sub V_\alpha^{(2)}:=V^{(2)}\cap V_\alpha$,
there exists $\beta\geq\alpha$ such that
$\conv_2(K)\sub V_\beta$.
\end{description}
Then
$U^{(2)}:=\phi^{-1}(V^{(2)})$ is an open
subset of~$U$, called a \emph{core} of~$\phi$.
If cores of well-filled charts
cover~$M$, then $M$ is said to \emph{admit
well-filled charts}.}
\end{definition}
On a first reading, the
reader may find the
notion of a well-filled chart
somewhat elusive.
Special cases of particular interest
(which are more concrete and
easier to understand)
are described in \cite[Examples~1.11 and 1.12]{HOM}.
See \cite[Theorem~1.13]{HOM}
for the following result.\index{homotopy group}
\begin{theorem}\label{wellfilled}
Let $M$ be a Hausdorff topological
space containing a
directed union $M_\infty:=\bigcup_{\alpha\in A}M_\alpha$
of Hausdorff topological spaces $M_\alpha$
as a dense subset, such that all inclusion maps
$M_\alpha\to M_\beta$ $($for $\alpha\leq \beta)$
and $M_\alpha\to M$ are continuous.
If $M$ admits well-filled
charts, then
\[
\pi_k(M,p)\;=\;\dl\, \pi_k(M_\alpha,p)\quad\mbox{for all $\,k\in \N_0$
and $p\in M_\infty$.}
\]
\end{theorem}
For a typical application, let $H$ be a Lie group,
$m\in \N$, $\cS(\R^m,\Lie(H))$ be the Schwartz space
of rapidly decreasing $\Lie(H)$-valued
smooth functions on~$\R^m$,
and $\cS(\R^m,H)$ be the corresponding Lie group,
as in \cite{BCR} (for special~$H$) and \cite{Wa2}.
Then $C^\infty_c(\R^m,H)=\bigcup_{n\in \N}C^\infty_{[{-n},n]^m}(\R^m,H)$
is dense in $\cS(\R^m,H)$, and $\cS(\R^m,H)$
admits well-filled charts~\cite[Example~8.4]{HOM}.
Using Theorem~\ref{wellfilled}
and approximation results from~\cite{NeG},
it is then easy to see that
\begin{eqnarray*}
\pi_k(\cS(\R^m,H)) \, &\isom & \, \dl\, \pi_k(C^\infty_{[{-n},n]^m}(\R^m,H))\;\isom \;
\pi_k(C^\infty_c(\R^m,H))\\
\,& \isom &\,
\pi_k(C_0(\R^m,H))\;\isom\; \pi_k(C(\bS_m,H)_*)\;\isom\;
\pi_{k+m}(H)
\end{eqnarray*}
(see \cite[Remark~8.6]{HOM}).
This had been conjectured in
\cite{BCR} and was open since 1981.
\section{Subgroups of ascending unions and related topics}\label{secsub}
We now discuss various results concerning subgroups
of ascending unions of Lie groups (notably for
direct limits
of finite-dimensional Lie groups).\\[3mm]
{\bf Non-existence of small subgroups.}
It is an open problem whether infinite-dimensional Lie groups
may contain small torsion subgroups~\cite[p.\,293]{NeS}.
For direct limits of finite-dimensional Lie groups,
the pathology could be ruled out by proving that they do not
contain small subgroups \cite[Theorem~A]{SMA}:\index{small subgroup}
\begin{theorem}
If $G_1\sub G_2\sub \cdots$
is a direct sequence of finite-dimensional Lie groups,
then the Lie group
$G=\dl\,G_n$\vspace{-.7mm}
does not have small subgroups.
\end{theorem}
\emph{Idea of proof.}
Given a compact identity neighbourhood
$C_1\sub G_1$ which does not contain non-trivial
subgroups of~$G_1$,
there exists a compact identity neighbourhood
$C_2\sub G_2$ with $C_1$ in its interior
relative $G_2$, which does not contain
non-trivial subgroups of $G_2$ (see \cite[Lemma~2.1]{SMA}).
Proceeding in this way, we find
a sequence $(C_n)_{n\in \N}$ of
compact identity neighbourhoods
$C_n\sub G_n$ not containing non-trivial subgroups,
such that $C_n\sub C_{n+1}^0$ for each~$n$.
Then $C:=\bigcup_{n\in \N}C_n$
is an identity neighbourhood in~$G$
and we may hope that~$C$ does not contain non-trivial
subgroups of~$G$. Unfortunately, this
is not true in general, as the example
$\R^{(\N)}=\bigcup_{n\in \N}\R^n=\bigcup_{n\in\N}C_n$
with $C_n:=[{-n},n]^n$ shows.
However, if the sets $C_n$
are chosen carefully (which requires much work),
then indeed $C$ will
not contain non-trivial
subgroups~\cite{SMA}.\vspace{2mm}\smartqed\qed

\noindent
We mention that an analogous result is available for certain
ascending unions
of infinite-dimensional Lie groups
$G_1\sub G_2\sub \cdots$ (see
\cite[Theorem~B]{SMA}).
To enable compactness arguments,
each $G_n$ has to be locally $k_\omega$
or each $G_n$ a Banach--Lie group
and the tangent map
$\Lie(\lambda_n)\colon \Lie(G_n)\to \Lie(G_{n+1})$
of the inclusion map $\lambda_n\colon
G_n\to G_{n+1}$ a compact operator.\footnote{Further technical
hypotheses need to be imposed,
which we suppress here.}\\[3mm]
{\bf Initial Lie subgroups.}
If $G$ is a Lie group and $H\sub G$ a subgroup, then
$H$ is called an \emph{initial Lie subgroup}\footnote{Some readers
may prefer to omit the second condition,
or allow $M$ to be a manifold with $C^k$-boundary, with corners
or (more generally) a $C^k$-manifold with rough
boundary (as introduced in \cite{GaN}).\index{initial Lie subgroup}
The following results carry over to these
varied situations (see \cite{OPE}).}
if it admits a Lie group structure
making the inclusion map $\iota \colon H\to G$
a smooth map,
such that $\Lie(\iota)$ is injective
and mappings from $C^k$-manifolds~$M$ to~$H$
are $C^k$ if and only if they are $C^k$ as
mappings to~$G$,
for each $k\in \N\cup\{\infty\}$.\\[2.5mm]
Answering an open problem from \cite{NeS}
in the negative, it was shown in~\cite{OPE}
that subgroups of
infinite-dimensional Lie groups
not be initial Lie subgroups.
In fact, one can take $G=\R^\N$ (with the product topology)
and $H=\ell^\infty$ (see \cite[Theorem~1.3]{OPE}).
For direct limits
of finite-dimensional Lie groups,
$G=\bigcup_{n\in \N}G_n$,
it was already shown in~\cite{FUN}
that every subgroup~$H\sub G$
admits a natural Lie group structure.
By \cite[Theorem~2.1]{OPE},
this Lie group structure makes~$H$
an initial Lie subgroup of~$G$
and thus the preceding pathology
does not occur for such
direct limit Lie groups~$G$.\vfill\pagebreak

\noindent
{\bf Continuous one-parameter groups and the
topology on {\boldmath$\Lie(G)$}.}
If $G=\bigcup_{n\in \N}G_n$
is a direct limit of finite-dimensional Lie groups,
then every continuous homomorphism
$(\R,+)\to G$ (i.e., each continuous one-parameter subgroup)
is a continuous homomorphism to some
$G_n$ (by compact regularity) and hence smooth.
It easily follows from this that
the natural map
\[
\theta\colon \Lie(G)\to \Hom_{\cts}(\R,G)\,,\quad
x\mto (t\mto \exp_G(tx))
\]
is a bijection onto the set $\Hom_{\cts}(\R,G)$
of continuous one-parameter subgroups of~$G$.
It was asked in \cite[Problem~VII.2]{NeS}
whether~$G$ is a \emph{topological group
with Lie algebra}\index{topological group
with Lie algebra} in the sense of \cite[Definition~2.11]{HaM}.
This holds if
$\theta$
is a homeomorphism onto
$\Hom_{\cts}(\R,G)$, equipped with the compact-open
topology
(which is not obvious because $\exp_G$ need
not be a local homeomorphism at~$0$).
As shown in \cite[Theorem~3.4]{OPE},
the latter property is always satisfied.
Thus $\Lie(G)$ is determined by the topological group
structure of~$G$. E.g., this implies
that every continuous homomorphism
from a locally exponential Lie group to~$G$
is smooth \cite[Proposition~3.7]{OPE}
(where a Lie group
is called \emph{locally exponential}\index{locally exponential Lie group}
if it has an exponential function
and the latter is a local
diffeomorphism at~$0$).\index{automatic smoothness}
It is an
open problem
whether
continuous homomorphisms between arbitrary Lie groups
are automatically smooth.
\begin{acknowledgement}
The author thanks K.-D. Bierstedt
and S.\,A. Wegner (Paderborn) for discussions
related to regularity properties
of (LF)-spaces, and J. Bonet (Valencia)
for comments which entered into Remark~\ref{bonet}.
K.-H. Neeb (Darmstadt) contributed useful comments
on an earlier version of the article.
The research was supported
by the German Research Foundation (DFG),
projects GL 357/5-1 and GL 357/7-1.
\end{acknowledgement}
\printindex

\begin{thebibliography}{9999MM}
%
%
\bibitem[ACM89]{ACM89}
M.\,C. Abbati, R. Cirelli, A. Mania and
P.\,W. Michor,
{\it The Lie group of automorphisms of
a principal bundle}, J. Geom.\ Phys.\ \textbf{6:2} (1989),
215--235.
%
%
\bibitem[AHK93]{AH93}
S.\,A. Albeverio,
R.\,J. H\o{}egh-Krohn, J.\,A. Marion,
D.\,H. Testard and B.\,S. Torr\'{e}sani,
``Noncommutative Distributions,''
Marcel Dekker, New York, 1993.
%
%
\bibitem[AK06]{AaK}
S.\,A. Albeverio and A. Kosyak,
{\it Quasiregular representations of the infinite-dimensional
nilpotent group},
J.\ Funct.\ Anal.\ \textbf{236:2} (2006), 634--681.
%
%
\bibitem[BBZ04]{BBZ04}
Y.\,A. Bahturin, A.\,A. Baranov and A.\,E. Zalesski,
{\it Simple Lie subalgebras of locally finite associative
algebras}, J. Algebra  {\bf 281:1} (2004), 225--246.
%
%
\bibitem[BB04]{BB4}
Y. Bahturin and G. Benkart,
{\it Some constructions in the theory of locally finite simple
Lie algebras},
J. Lie Theory \textbf{14:1} (2004), 243--270. 
%
%
\bibitem[Ba97]{Ban}
A. Banyaga, ``The Structure of Classical
Diffeomorphism Groups,''
Kluwer, Dordrecht, 1997.
%
%
\bibitem[Bi88]{Bie}
K.-D. Bierstedt,
{\it An introduction to locally convex inductive limits},
pp.\,35--133
in: H.\ Hogbe-Nlend (ed.), ``Functional Analysis and its Applications,''
World Scientific, Singapore, 1988.
%
%
\bibitem[BMS82]{BMS}
K.-D. Bierstedt, R. Meise
and W.\,H. Summers,
{\it A projective description of weighted inductive
limits}, Trans.\ Amer.\ Math.\ Soc.\ {\bf 272:1}
(1982), 107--160.
%
%
\bibitem[BS71]{BaS}
J. Bochnak and J. Siciak, 
Analytic functions in topological vector spaces.
Studia Math.\ \textbf{39} (1971), 77--112.
%
%
\bibitem[BCR81]{BCR}
H. Boseck, G. Czichowski and K.-P. Rudolph,
``Analysis on Topological Groups -- General Lie Theory,''
Teubner, Leipzig, 1981.
%
%
\bibitem[Bo87]{BTV}
N. Bourbaki, ``Topological Vector Spaces,
Chapters 1--5,'' Springer, Berlin, 1987.
%
%
\bibitem[Ch85]{Ch85}
S.\,B. Chae, ``Holomorphy and Calculus in Normed
Spaces,'' Marcel Dekker, New York, 1985.
%
%
\bibitem[Da08]{Dah}
R. Dahmen,
{\it Complex analytic mappings on}
(\hspace*{-.2mm}{\it LB}\hspace*{.3mm})-{\it spaces
and applications in infinite-dimensional
Lie theory},
in preparation.
%
%
\bibitem[DG08]{REG}
R. Dahmen
and H. Gl\"{o}ckner,
{\it Regularity in Milnor's sense
for direct limits of infinite-dimensional
Lie groups},
in preparation.
%
%
\bibitem[DW97]{DaW}
S. Dierolf and J. Wengenroth,
{\it Inductive limits of topological algebras},
Linear Topol.\ Spaces Complex Anal.\ \textbf{3} (1997),
45--49.
%
%
\bibitem[DP04]{DP04}
I. Dimitrov and I. Penkov,
{\it Borel subalgebras of $\gl(\infty)$},
Resenhas {\bf 6}:2--3 (2004), 153--163. 
%
%
\bibitem[DPW02]{DPW}
I. Dimitrov, I. Penkov and J.\,A. Wolf,
{\it A Bott-Borel-Weil theory for direct limits of algebraic
groups}, Amer.\ J. Math.\ \textbf{124:5} (2002), 955--998.
%
%
\bibitem[Du64]{Dud} R.\,M. Dudley,
{\it On sequential convergence},
Trans.\ Amer.\ Math.\ Soc.\ \textbf{112} (1964), 483--507.
%
%
\bibitem[Dv02]{Dv02}
A. Dvorsky, {\it Direct limits of low-rank representations of
classical groups}, Comm.\ Algebra \textbf{30:12} (2002),
6011--6022.
%
%
\bibitem[Ed99]{Eda}
T. Edamatsu, {\it On the bamboo-shoot topology of certain
inductive limits of topological groups},
J. Math.\ Kyoto Univ.\ \textbf{39:4}, 715--724.
%
%
\bibitem[Fl71]{Fl71}
K. Floret,
{\it Lokalkonvexe Sequenzen mit kompakten Abbildungen},
J. Reine Angew.\ Math.\ \textbf{247} (1971), 155--195.
%
%
\bibitem[Fl79]{Fl79}
K. Floret,
{\it On bounded sets in inductive limits of normed spaces}, 
Proc.\ Amer.\ Math.\ Soc.\  \textbf{75:2} (1979),
221--225. 
%
%
\bibitem[Fl80]{Fl80}
K. Floret,
{\it Some aspects of the theory of locally convex
inductive limits}, pp.\,205--237 in:
K.-D. Bierstedt and B. Fuchssteiner (eds.),
``Functional Analysis: Surveys and Recent Results II,''
North-Holland Math.\ Stud.\ \textbf{38},
North-Holland, Amsterdam, 1980.
%
%
\bibitem[Gl02a]{RES}
H. Gl\"{o}ckner,
{\it Infinite-dimensional Lie groups without
completeness restrictions},
pp.\,43--59
in: A. Strasburger et al.\ (eds.),
``Geometry and Analysis on Finite- and Infinite-Dimensional
Lie Groups,''
Banach Center Publ.\ \textbf{55}, Warsaw, 2002.
%
%
\bibitem[Gl02b]{GCX}
H. Gl\"{o}ckner,
{\it Lie group structures on quotient groups
and universal complexifications for infinite-dimensional
Lie groups},
J. Funct.\ Anal.\ \textbf{194:2} (2002), 347--409.
%
%
\bibitem[Gl02c]{ALG}
H. Gl\"{o}ckner,
{\it Algebras whose groups of units are Lie groups},
Studia Math.\ {\bf 153:2} (2002), 147--177.
%
%
\bibitem[Gl02d]{DIF}
H. Gl\"{o}ckner,
{\it Patched locally convex spaces,
almost local mappings
and diffeomorphism groups of non-compact
manifolds}, manuscript, 2002.
%
%
\bibitem[Gl03a]{DIR}
H. Gl\"{o}ckner,
{\it Direct limit Lie groups and manifolds},
J. Math.\ Kyoto Univ.\ \textbf{43:1} (2003), 1--26.
%
%
\bibitem[Gl03b]{MEA}
H. Gl\"{o}ckner,
{\it Lie groups of measurable mappings},
Canad.\ J. Math.\ \textbf{55:5} (2003), 969--999.
%
%
\bibitem[Gl04a]{HOL}
H. Gl\"{o}ckner,
{\it Lie groups of germs of analytic mappings},
pp.\,1--16
in: T. Wurzbacher (ed.),
``Infinite-Dimensional Groups and Manifolds,''
IRMA Lecture Notes in Math.\ and Theor.\ Physics,
de Gruyter, 2004.
%
%
\bibitem[Gl04b]{ZOO}
H. Gl\"{o}ckner, {\it Lie groups
over non-discrete topological fields},
preprint, arXiv:math/0408008.
%
%
\bibitem[Gl05a]{DFR}
H. Gl\"{o}ckner,
{\it $\Diff(\R^n)$ as a Milnor-Lie group},
Math.\ Nachr.\ \textbf{278:9} (2005),
1025--1032.
%
%
\bibitem[Gl05b]{FUN}
H. Gl\"{o}ckner,
{\it Fundamentals of direct limit Lie theory},
Compos.\ Math.\ \textbf{141:6} (2005), 1551--1577.
%
%
\bibitem[Gl07a]{SMA}
H. Gl\"{o}ckner,
{\it Direct limit groups do not have small subgroups},
Topol.\ Appl.\ \textbf{154:6} (2007), 1126--1133.
%
%
\bibitem[Gl07b]{COM}
H. Gl\"{o}ckner,
{\it Direct limits of infinite-dimensional Lie groups
compared to direct limits in related categories},
J. Funct.\ Anal.\ \textbf{245:1} (2007), 19--61.
%
%
\bibitem[Gl08a]{OPE}
H. Gl\"{o}ckner,
{\it Solutions to questions in Neeb's recent survey on
infinite-dimensional Lie groups},
preprint, arXiv:0801.1919v1.
%
%
\bibitem[Gl08b]{HOM}
H. Gl\"{o}ckner,
{\it Homotopy groups of ascending unions
of infinite-dimensional manifolds},
in preparation.
%
%
\bibitem[GGH07]{GGH}
H. Gl\"{o}ckner, R. Gramlich and T. Hartnick,
{\it Final group topologies,
Phan systems and Pontryagin duality},
preprint, arXiv:math.GR/0603537.
%
%
\bibitem[GN08]{GaN}
H. Gl\"{o}ckner and K.-H. Neeb,
``Infinite-Dimensional Lie Groups, Vol~1,''
book in preparation.
%
%
\bibitem[Go03]{Gol}
G.\,A. Goldin,
{\it Lectures on diffeomorphism groups in
quantum physics},
pp.\,3--93
in:
``Contemporary Problems in Mathematical Physics,''
Proceedings of the 3rd Intern.\ Workshop
(Cotonou, 2003), World Scientific, 2004
%
%
\bibitem[Ha71]{Han}
V.\,L. Hansen, {\it Some theorems on direct limits
of expanding systems of manifolds},
Math.\ Scand.\ \textbf{29} (1971), 5--36.
%
%
\bibitem[Hi93]{Hi93}
T. Hirai,
{\it Irreducible unitary representations of the group of
diffeomorphisms of a noncompact manifold},
J. Math.\ Kyoto Univ.\ \textbf{33:3} (1993), 827--864.
%
%
\bibitem[HST01]{HSTH}
T. Hirai, H. Shimomura, N. Tatsuuma and E. Hirai,
{\it Inductive limits of topologies, their direct product,
and problems related to algebraic structures},
J. Math.\ Kyoto Univ.\ \textbf{41:3} (2001), 475--505.
%
%
\bibitem[HM07]{HaM}
K.\,H. Hofmann and S.\,A. Morris,
``The Lie Theory of Connected Pro-Lie Groups,''
EMS Tracts in Mathematics \textbf{2},
European Mathematical Society, Zurich, 2007.
%
%
\bibitem[Is76]{Ism}
R.\,S. Ismagilov,
{\it Unitary representations of the group $C^\infty_0(X,G)$,
$G=\SU_2$}, 
Mat.\ Sb.\ (N.S.)  \textbf{100(142):1} (1976),
117--131, 166. 
%
%
\bibitem[KR01]{KaR}
N. Kamran and T. Robart,
{\it A manifold structure for analytic isotropy
Lie pseudogroups of infinite type},
J. Lie Theory \textbf{11:1} (2001), 57--80.
%
%
\bibitem[Ke74]{Kel}
H.\,H. Keller,
``Differential Calculus in Locally Convex Spaces,''
Springer, Berlin, 1974.
%
%
\bibitem[KOV04]{KOV}
S. Kerov, G.\,I. Ol'shanski\u{\i} and
A. Vershik,
{\it Harmonic analysis on the infinite symmetric group},
Invent.\ Math.\ \textbf{158:3} (2004), 551--642.
%
%
\bibitem[Ki81]{Ki81}
A.\,A. Kirillov,
{\it Unitary representations of the group
of diffeomorphisms and of some of its subgroups},
Sel.\ Math.\ Sov.\ \textbf{1} (1981), 351--372.
%
%
\bibitem[Ko69]{Ko69}
J. K\"{o}hn,
{\it Induktive Limiten nicht lokal-konvexer topologischer
Vektorr\"{a}ume}, 
Math.\ Ann.\ \textbf{181} (1969), 269--278.
%
%
\bibitem[KS77]{KS77}
V.\,I. Kolomycev and Y.\,S. Samo\u{\i}lenko,
{\it Irreducible representations of inductive limits of groups},
Ukr.\ Mat.\ J. \textbf{29} (1977), 526--531, 565
%
%
\bibitem[KM97]{KaM}
A. Kriegl and P.\,W. Michor,
``The Convenient Setting of
Global Ana\-lysis,'' Amer.\ Math.\ Soc.,
Providence, 1997.
%
%
\bibitem[Ku06]{Ku06}
K. K\"{u}hn,
{\it Direct limits of diagonal chains of type O, U, and Sp,
and their homotopy groups},
Comm.\ Algebra \textbf{34:1} (2006), 75--87.
%
%
\bibitem[Le94]{Lei}
F. Leitenberger,
{\it Unitary representations and coadjoint orbits
for a group of germs of real analytic diffeomorphisms},
Math.\ Nachr.\ \textbf{169} (1994),
185--205.
%
%
\bibitem[LS00]{LS00}
A.\,K. Leonov and S.\,A. Shkarin,
{\it Integral completeness of locally convex spaces},
Russ.\ J. Math.\ Phys.\ \textbf{7:4} (2000), 402--412.
%
%
\bibitem[Ls82]{Les}
J. Leslie,
{\it On the group of real analytic diffeomorphisms of a
compact real analytic manifold},
Trans.\ Amer.\ Math.\ Soc.\ \textbf{274:2} (1982),
651--669.
%
%
\bibitem[Ls85]{Ls85}
J. Leslie,
{\it Some finite-codimensional Lie subgroups of $\Diff^\omega(M)$},
pp.\,359--372 in: G.\,M. Rassias and T.\,M. Rassias (eds.),
``Differential Geometry, Calculus of Variations, and their Applications,''  
Lecture Notes in Pure and Appl.\ Math.\ {\bf 100},
Marcel Dekker, New York, 1985. 
%
%
\bibitem[Mk63]{Mk63}
B.\,M. Makarov,
{\it Some pathological properties of inductive limits of $B$-spaces},
Uspehi Mat.\ Nauk  \textbf{18:3} (1963), 171--178. 
%
%
\bibitem[Mc80]{Mic}
P.\,W. Michor,
``Manifolds of Differentiable Mappings,''
Shiva Publishing, Nantwich, 1980.
%
%
\bibitem[Mr82]{Mr82}
J. Milnor,
{\it On infinite dimensional Lie groups},
Preprint, Institute of\linebreak
Advanced Study, Princeton,
1982.
%
%
\bibitem[Mr84]{Mr84}
J. Milnor,
{\it Remarks on infinite dimensional Lie groups},
pp.\,1007--1057 in:
B. DeWitt and R. Stora (eds.),
``Relativity, groups and topology II,''
North-Holland, Amsterdam, 1984.
%
%
\bibitem[Mi66]{Miu}
A.\,A. Miljutin,
{\it Isomorphism of the spaces of continuous functions over compact sets
of the cardinality of the continuum},
Teor.\ Funkci\u{\i} Funkcional.\ Anal.\ i Prilo\v{z}en.\ Vyp.\  
\textbf{2} (1966), 150--156.
%
%
\bibitem[Mj79]{Mj79}
J. Mujica,
{\it Spaces of germs of holomorphic functions},
pp.\,1--41 in:
G.-C. Rota (ed.),
``Studies in Analysis,''  
Adv.\ in Math.\ Suppl.\ Stud.\ \textbf{4},
Academic Press, New York, 1979. 
%
%
\bibitem[Mj83]{Muj}
J. Mujica,
\emph{Spaces of continuous functions with values in an inductive
limit}, pp.\,359--367
in: G.\,I. Zapata (ed.),
``Functional Analysis, Holomorphy, and Approximation theory,''
Lect.\ Notes Pure Appl.\ Math.\ \textbf{83},
Marcel Dekker, New York, 1983.
%
%
\bibitem[NRW91]{NRW1}
L. Natarajan, E. Rodr\'{\i}guez-Carrington and J.\,A.
Wolf,
{\it Differentiable structure for direct limit groups},
Letters in Math.\ Phys.\ \textbf{23:2}, (1991)
99--109.
%
%
\bibitem[NRW93]{NRW2}
L. Natarajan, E. Rodr\'{\i}guez-Carrington and
J.\,A. Wolf,
{\it Locally convex Lie groups},
Nova J.\ Alg.\ Geom.\ \textbf{2:1} (1993),
59--87.
%
%
\bibitem[NRW01]{NRW4}
L. Natarajan, E. Rodr\'{\i}guez-Carrington and
J.\,A. Wolf,
{\it The Bott-Borel-Weil theorem for direct limit groups},
Trans.\ Amer.\ Math.\ Soc.\ \textbf{353:11} (2001),
4583--4622.
%
%
\bibitem[Ne98]{Ne98}
K.-H. Neeb,
{\it Holomorphic highest weight representations of infinite-dimensional
complex classical groups},
J. Reine Angew.\ Math.\ \textbf{497} (1998),
171--222.
%
%
\bibitem[Ne00]{Ne00}
K.-H. Neeb,
{\it Integrable roots in split graded Lie algebras},
J. Algebra  \textbf{225:2} (2000), 534--580.
%
%
\bibitem[Ne02a]{NeH}
K.-H. Neeb,
{\it Classical Hilbert-Lie groups, their extensions
and their homotopy groups},
pp.\,87--151 in: A. Strasburger et al.\ (eds.)
``Geometry and Analysis on Finite- and Infinite-Dimensional
Lie Groups,'' Banach Center Publ.\ \textbf{55},
Warsaw, 2002.
%
%
\bibitem[Ne02b]{NeC}
K.-H. Neeb,
{\it Central extensions of infinite-dimensional
Lie groups}, Ann.\ Inst.\ Fourier (Grenoble)
\textbf{52:5} (2002), 1365--1442.
%
%
\bibitem[Ne04a]{NeG}
K.-H. Neeb,
{\it Current groups for non-compact manifolds and
their central extensions},
pp.\,109--183 in:
T. Wurzbacher (ed.)
``Infinite-Dimensional Groups and Manifolds,''
IRMA Lecture Notes in Math.\ and Theor.\ Physics,
de Gruyter, 2004.
%
%
\bibitem[Ne04b]{NeA}
K.-H. Neeb,
{\it Abelian extensions
of infinite-dimensional Lie groups},
Travaux Math.\ \textbf{15} (2004),
69--194.
%
%
\bibitem[Ne04c]{NeO}
K.-H. Neeb,
{\it Infinite-dimensional groups and their representations},
pp.\,213--328 in: J.\,C. Jantzen and K.-H. Neeb (eds.),
``Lie Theory,''
Birkh\"{a}user, Boston, 2004.
%
%
\bibitem[Ne06]{NeS}
K.-H. Neeb,
{\it Towards a Lie theory of locally convex groups},
Jpn.\ J. Math.\ \textbf{1:2} (2006), 291--468.
%
%
\bibitem[Ne07]{NeN}
K.-H. Neeb,
{\it Non-abelian extensions of infinite-dimensional
Lie groups}, Ann.\ Inst.\ Fourier
(Grenoble) \textbf{57:1} (2007),
209--271.
%
%
\bibitem[NS01]{NS01}
K.-H. Neeb and N. Stumme,
{\it The classification of locally finite split simple Lie algebras},
J. Reine Angew.\ Math.\textbf{533} (2001), 25--53. 
%
%
\bibitem[NW07]{NaW}
K.-H. Neeb and F. Wagemann,
{\it Lie group structures on groups of smooth and
holomorphic maps on non-compact manifolds},
to appear in Geom.\ Dedicata
(cf.\ arXiv:math/0703460).
%
%
\bibitem[Ol83]{Ol83}
G.\,I. Ol'shanski\u{\i},
{\it Unitary representations of infinite-dimensional pairs
$(G,K)$ and the formalism of R. Howe},
Dokl.\ Akad.\ Nauk SSSR \textbf{269:1}
(1983), 33--36.
%
%
\bibitem[Ol03]{Ol03}
G.\,I. Ol'shanski\u{\i},
{\it The problem of harmonic analysis on the infinite-dimensional
unitary group}, J. Funct.\ Anal.\ \textbf{205:2}
(2003), 464--524.
%
%
\bibitem[Om97]{Om97}
H. Omori,
``Infinite-Dimensional Lie Groups,''
Transl.\ Math.\ Monographs \textbf{158},
Amer.\ Math.\ Soc., Providence, 1997.
%
%
\bibitem[Pa65]{Pa1}
R.\,S. Palais,
{\it On the homotopy type
of certain groups of operators},
Topology \textbf{3} (1965),
271--279.
%
%
\bibitem[Pa66]{Pa2}
R.\,S. Palais,
{\it Homotopy theory
of infinite-dimensional manifolds},
Topology \textbf{5} (1966),
1--16.
%
%
\bibitem[PS03]{PS03}
I. Penkov and H. Strade,
{\it Locally finite Lie algebras with root decomposition},
Arch.\ Math.\ (Basel) \textbf{80:5} (2003), 478--485. 
%
%
\bibitem[Pk08]{Pi08}
D. Pickrell
{\it Heat kernel measures and critical limits},
also published in this volume.
%
%
\bibitem[Pi77]{Pis}
D. Pisanelli,
{\it An example of an infinite Lie group},
Proc.\ Amer.\ Math.\ Soc.\
\textbf{62:1} (1977),
156--160.
%
%
\bibitem[Ra07]{Ra07}
M. Rabaoui,
{\it A Bochner type theorem for inductive limits of
Gelfand pairs},
preprint, arXiv:math/0701501v1.
%
%
\bibitem[Sm78]{Sme}
J. Schmets,
{\it Spaces of vector-valued continuous functions},
pp.\,368--377 in: R.\,M. Aron and S. Dineen (eds.),
``Vector Space Measures and Applications II,''
Lecture Notes in Math.\ \textbf{644},
Springer, Berlin, 1978.
%
%
\bibitem[Sm83]{Sm83}
J. Schmets,
``Spaces of Vector-Valued Continuous Functions,''
Lecture Notes in Math.\ \textbf{1003},
Springer, Berlin, 1983.
%
%
\bibitem[Se55]{Si55} 
J. Sebasti\~{a}o e Silva,
{\it Su certe classi di spazi localmente convessi importanti per le
applicazioni},
Rend.\ Mat.\ Appl.\ \textbf{14} (1955), 388--410.
%
%
\bibitem[Sh01]{Sh01}
H. Shimomura,
{\it Quasi-invariant measures on the group of
diffeomorphisms and smooth vectors
of unitary representations},
J. Funct.\ Anal.\ \textbf{187:2} (2001),
406--441.
%
%
\bibitem[Sh05]{Sh05}
H. Shimomura,
{\it Irreducible decompositions of unitary representations
of infinite-dimensional groups},
Math.\ Z. \textbf{251:3} (2005),
575--587.
%
%
\bibitem[Sr59]{Sr59}
T. Shirai,
{\it Sur les topologies des espaces de L. Schwartz}, 
Proc.\ Japan Acad.\  \textbf{35} (1959), 31--36.
%
%
\bibitem[Sv69]{Sm69}
O.\,G. Smolyanov,
{\it Almost closed linear subspaces of strict inductive limits of sequences
of Fr\'{e}chet spaces},
Mat.\ Sb.\ (N.S.)  \textbf{80 (122)} (1969), 513--520.
%
%
\bibitem[St99]{St99}
N. Stumme,
{\it The structure of locally finite split Lie algebras},
J. Algebra {\bf 220:2} (1999), 664--693.
%
%
\bibitem[TSH98]{TSH}
N. Tatsuuma, H. Shimomura and T. Hirai,
{\it On group topologies and unitary representations of inductive
limits of topological groups and the case of the group of
diffeomorphisms},
J. Math.\ Kyoto Univ.\ \textbf{38:3} (1998),
551--578.
%
%
\bibitem[Th64a]{Th64a}
E. Thoma,
{\it \"{U}ber unit\"{a}re Darstellungen abz\"{a}hlbarer, diskreter
Gruppen},
Math.\ Ann.\ \textbf{153} (1964), 111--138.
%
%
\bibitem[Th64b]{Th64b}
E. Thoma,
{\it Die unzerlegbaren, positiv-definiten Klassenfunktionen der abz\"{a}hlbar
unendlichen, symmetrischen Gruppe},
Math.\ Z. \textbf{85} (1964), 40--61.
%
%
\bibitem[VGG75]{VGG75}
A.\,M. Ver\v{s}ik, I.\,M. Gel'fand and M.\,I. Graev,
{\it Representations of the group of diffeomorphisms},
Uspehi Mat.\ Nauk\ \textbf{30:6} (1975), 1--50.
%
%
\bibitem[Vo76]{Voi}
D. Voiculescu,
{\it Repr\'{e}sentations factorielles de type $II_1$ de $\U(\infty)$},
J. Math.\ Pures Appl.\ \textbf{55:1} (1976),
1--20.
%
%
\bibitem[Wa08]{Wa2}
B. Walter,
{\it Weighted diffeomorphism groups of Banach spaces
and weighted mapping groups},
in preparation.
%
%
\bibitem[We03]{Wen}
J. Wengenroth,
``Derived Functors in Functional Analysis,''
Lecture Notes in Math.\ \textbf{1810},
Springer, Berlin, 2003.
%
%
\bibitem[Wk07]{Wk07}
C.\ Wockel,
{\it Lie group structures on symmetry groups of principal bundles},
J. Funct.\ Anal.\ \textbf{251:1} (2007),
254--288.
%
%
\bibitem[Wo05]{Wo05}
J.\,A. Wolf,
{\it Principal series representations of direct limit
groups},
Compos.\ Math.\ {\bf 141:6} (2005),
1504--1530.
%
%
\bibitem[Wo08]{Wo08}
J.\,A. Wolf,
{\it
Infinite dimensional multiplicity free spaces} I:
{\it limits of compact commutative spaces},
also published in this volume.
%
%
\bibitem[Ya98]{Yam}
A. Yamasaki,
{\it Inductive limits of general linear groups},
J. Math.\ Kyoto Univ.\ {\bf 38:4} (1998),
769--779.
%
%
\end{thebibliography}
\end{document}